# FIXED-DOMAIN ASYMPTOTIC PROPERTIES OF TAPERED MAXIMUM LIKELIHOOD ESTIMATORS

By Juan Du, Hao Zhang[1] and V. S. Mandrekar

*Michigan State University, Purdue University and Michigan State University*

When the spatial sample size is extremely large, which occurs in many environmental and ecological studies, operations on the large covariance matrix are a numerical challenge. Covariance tapering is a technique to alleviate the numerical challenges. Under the assumption that data are collected along a line in a bounded region, we investigate how the tapering affects the asymptotic efficiency of the maximum likelihood estimator (MLE) for the microergodic parameter in the Matérn covariance function by establishing the fixed-domain asymptotic distribution of the exact MLE and that of the tapered MLE. Our results imply that, under some conditions on the taper, the tapered MLE is asymptotically as efficient as the true MLE for the microergodic parameter in the Matérn model.

**1. Introduction.** With the advancement of technology, large amounts of data are routinely collected over space and/or time in many studies in environmental monitoring, climatology, hydrology and other fields. The large amounts of correlated data present a great challenge to the statistical analysis and may render some traditional statistical approaches impractical. For example, in the maximum likelihood or Bayesian inference, the inverse of an $n \times n$ covariance matrix is involved, where the sample size $n$ may be in hundreds of thousands or even larger. Inverting the large covariance matrix repeatedly is a great computational burden if not impractical, and some approximation to the likelihood is necessary.

Covariance tapering is one of the approaches to approximating the covariance matrix and, therefore, the likelihood. Let the second order stationary Gaussian process $X(\mathbf{t}), \mathbf{t} \in \mathbb{R}^d$ have mean 0 and an isotropic covariance

Received June 2008; revised October 2008.
[1]Supported by NSF Grants DMS-04-05782 and DMS-07-06835.
*[AMS 2000 subject classifications.](https://www.ams.org/msc/)* Primary 62M20; secondary 62G20, 60G15.

function $K(h;\theta,\sigma^2)$, where $\sigma^2$ is the variance of the process and $\theta$ is the parameter that controls how fast the covariance function decays. Given $n$ observations $\mathbf{X}_n = (X(\mathbf{t}_1),\ldots,X(\mathbf{t}_n))'$, the log-likelihood is

$$(1.1) \quad l_n(\theta,\sigma^2) = -\frac{n}{2}\log 2\pi - \frac{1}{2}\log[\det \mathbf{V}_n(\theta,\sigma^2)] - \frac{1}{2}\mathbf{X}_n'[\mathbf{V}_n(\theta,\sigma^2)]^{-1}\mathbf{X}_n,$$

where $\mathbf{V}_n(\theta,\sigma^2)$ denotes the covariance matrix of $\mathbf{X}_n$. The idea of tapering is to keep the covariances approximately unchanged at small distance lags and to reduce the covariances to zero at large distances. To implement the idea, let $K_{\text{tap}}$ be an isotropic correlation function of compact support; that is, $K_{\text{tap}}(h) = 0$ if $h \geq \gamma$ for some $\gamma > 0$. Then, the tapered covariance function $\tilde{K}$ is the product of $K$ and $K_{\text{tap}}$,

$$(1.2) \quad \tilde{K}(h;\theta,\sigma^2) = K(h;\theta,\sigma^2)K_{\text{tap}}(h),$$

and the tapered covariance matrix is a Hadamard product $\tilde{\mathbf{V}}_n = \mathbf{V}_n(\theta,\sigma^2) \circ \mathbf{T}_n$, where $\mathbf{T}_n$ has the $(i,j)$th element as $K_{\text{tap}}(\|\mathbf{t}_i - \mathbf{t}_j\|)$. The tapered covariance matrix has a high proportion of zero elements and is, therefore, a sparse matrix. Inverting a sparse matrix is much more efficient computationally than inverting a regular matrix of the same dimension [see, e.g., Pissanetzky (1984), Gilbert, Moler and Schreiber (1992) and Davis (2006)]. One would use the tapered covariance function $\tilde{K}$ for spatial interpolation and estimation as if it was the correct covariance function. For example, the tapered maximum likelihood estimator maximizes the corresponding log-likelihood

$$(1.3) \quad l_{n,\text{tap}}(\theta,\sigma^2) = -\frac{n}{2}\log 2\pi - \frac{1}{2}\log[\det \tilde{\mathbf{V}}_n] - \frac{1}{2}\mathbf{X}_n'\tilde{\mathbf{V}}_n^{-1}\mathbf{X}_n.$$

Intuitively, if the taper is sufficiently close to 1 in the neighborhood of the origin, the tapering would not change the behavior of the covariance function near the origin. It has long been known that the behavior of the covariance function near the origin is most important to spatial interpolation. Stein (1988, 1990a, 1990b, 1999a, 1999b) has established rigorous fixed-domain asymptotic theory for spatial interpolation. Applying the general fixed-domain asymptotic theory, Furrer, Genton and Nychka (2006) showed that appropriate tapering does not affect the fixed-domain asymptotic mean square error of prediction for Matérn model.

Kaufman, Schervish and Nychka (2008) showed that the parameter in the Matérn covariance function, which is consistently estimable under the fixed-domain asymptotic framework, can be estimated consistently by the tapered MLE with $\theta$ fixed. However, it is unknown if the covariance tapering results in any loss of asymptotic efficiency.

The main objective of this paper is to establish the asymptotic properties, and particularly the asymptotic distribution of tapered MLE under the fixed-domain asymptotic framework. We now make a few remarks about why we adopt the fixed-domain asymptotic framework. When the spatial domain is fixed and bounded, more sample data can be obtained by sampling the domain increasingly densely. This results in the fixed-domain asymptotic framework. It is known that not all parameters in the covariance function are consistently estimable [e.g., Zhang (2004)] under the fixed-domain asymptotic framework. Zhang and Zimmerman (2005) argued that MLEs of the microergodic parameters are generally consistent but those of the nonmicroergodic parameters in general converge in distribution to a nondegenerate distribution. We refer readers to Stein (1999b), page 163, for the definition of microergodic parameters. In addition, Stein has established asymptotic results that show only the microergodic parameters affect the asymptotic mean square error under the fixed-domain asymptotic framework.

However, there is another asymptotic framework, where more data are sampled by increasing the spatial domain. This is the increasing domain asymptotic framework. Under mild regularity conditions, MLEs for all parameters are consistent and asymptotically normal [see Mardia and Marshall (1984)]. Therefore, asymptotic results are quite different under the two asymptotic frameworks. Given a finite sample, one has to know which asymptotic framework is more appropriate in order to apply any asymptotic results. Zhang and Zimmerman (2005) provided some guideline on this through both theoretical and numerical studies. Their results show that, for the exponential covariance function, the fixed-domain asymptotic distribution approximates the finite sample distribution at least as well as the increasing domain asymptotic distribution does. More specifically, for microergodic parameters, approximations corresponding to the two frameworks perform about equally well. For the nonmicroergodic parameters, the fixed-domain asymptotic approximation is preferable. In light of these results, we adopt the fixed-domain asymptotic framework in this work.

Fixed-domain asymptotic results for estimation are difficult to derive in general and there are only few results in literature [see Stein (1990c), Ying (1991, 1993), Chen, Simpson and Ying (2000), Zhang (2004), Loh (2005) and Kaufman, Schervish and Nychka (2008)]. Existing asymptotic distributions have been established only for specific models such as the exponential model for covariance functions [see Ying (1991, 1993) and Chen, Simpson and Ying (2000)] and a particular Matérn model with the smoothness parameter $\nu = 1.5$ [see Loh (2005)]. For the general Matérn model, the fixed-domain asymptotic distribution is not available even when data are observed along a line. In order to evaluate the efficiency of the tapered MLE, we establish the fixed-domain asymptotic distribution of MLE for the microergodic parameter in the general Matérn model [Theorem 5(i)] under the assumption that data are collected along a line. This result is of interest in its own right, outside the context of tapering.

It is even more difficult to study asymptotic properties of tapered MLE. Indeed, we are not aware of any fixed-domain asymptotic distribution established for tapered MLE. For this reason, we will start with a simple model, the Ornstein–Uhlenbeck process along a line, which is a stationary Gaussian process with zero mean and an exponential covariance function, and has Markovian properties. Due to the Markovian properties, the inverse of the covariance matrix can be given in closed form and is a band matrix. Therefore, for this model, it is not necessary to approximate the likelihood function. However, this simple model serves as a starting point in the study of covariance tapering and provides insight into the more general settings, which we will study subsequently.

Although spatial data are usually collected over a spatial region, there are situations when data are collected along lines. One example is the International H2O project, where measurements of meteorological data were collected by surface stations and aircraft along three flight paths that are along straight lines and transect the varied environmental conditions of the southern Great Plains [see Weckworth et al. (2004), LeMone et al. (2007) and Stassberg et al. (2008)]. Ecological data are sometimes collected along line transects as well.

The main results for the Ornstein–Uhlenbeck process are presented in Section 2. For the microergodic parameter in the Ornstein–Uhlenbeck process, we establish the asymptotic distribution of tapered MLE. In Section 3, we present the main results for a Gaussian stationary process having a Matérn covariogram. We put all proofs in Appendices A and B.

**2. Exponential model.** We assume the underlying process $X(t), t \in [0, 1]$ is Gaussian that has a mean 0 and an isotropic exponential covariogram $K(h) = \sigma^2 \exp(-\theta h)$. Such a process is known as the Ornstein–Uhlenbeck process, which has a Markovian property that will be exploited in our proof.

The exponential isotropic covariance function is one of the most commonly used models for spatial data analysis. It follows from Ying (1991) and Zhang (2004) that both $\sigma^2$ and $\theta$ are not consistently estimable under the fixed-domain asymptotic framework, but the product $\sigma^2 \theta$ is. Applying the fixed-domain asymptotic theory for spatial interpolation, Zhang (2004) showed that it is only this product, and not the individual parameters $\sigma^2$ and $\theta$, that asymptotically affects the interpolation. Therefore, it is important to estimate this product well. In this section, we establish the asymptotic properties of the tapered MLE of this product. For simplicity of argument, we will maximize the likelihood function over $(\theta, \sigma^2) \in J = [a, b] \times [w, v]$ for some constants $0 < a \leq b$ and $0 < w \leq v$ and do not require that $J$ contains the true parameter value $(\theta_0, \sigma_0^2)$. However, we do assume that $\theta_0 \sigma_0^2 \in \{\theta \sigma^2, (\theta, \sigma^2) \in J\}$; that is, there exists a pair $(\theta, \sigma^2)$ in $J$ such that $\theta \sigma^2 = \theta_0 \sigma_0^2$.

The following two assumptions are made throughout this section:

(A1) The process is observed at points $t_{k,n} \in [0,1], k = 1, \ldots, n$, with $0 \leq t_{1,n} < t_{2,n} < \cdots < t_{n,n} \leq 1$, and suppose that $n\Delta_{k,n}$ is bounded away from 0 and $\infty$, where $\Delta_{k,n} = t_{k,n} - t_{k-1,n}, k = 2, \ldots, n$. We also assume that $t_{n,n} \to 1$ and $t_{1,n} \to 0$ as $n \to \infty$.

(A2) $K_{\text{tap}}(h; \gamma)$ is an isotropic correlation function such that $K_{\text{tap}}(h; \gamma) = 0$ if $h \geq \gamma$, where $\gamma \in (0,1)$ is a constant. Moreover, $K_{\text{tap}}(h; \gamma)$ has a bounded second derivative in $h \in (0,1)$ and $K'_{\text{tap}}(h; \gamma) = ch + o(h)$ as $h \to 0+$ for some constant $c$.

A taper can be any correlation function with compact support, and such correlation functions have been studied in literature [see Wu (1995), Wendland (1995, 1998) and Gneiting (1999, 2002)]. We believe that a large number of compactly supported correlation functions satisfy assumption (A2). Particularly, a Wendland taper is a truncated polynomial and, therefore, satisfies (A2) if the degree of the polynomial is greater than 3.

We also note that the assumption in (A2) that $K_{\text{tap}}$ has a bounded second derivative in $h \in (0,1)$ can be weakened, so that $d^2 K_{\text{tap}}/dh^2$ exists at any $h \in (0, \gamma)$ as long as the first derivative exists everywhere in $(0,1)$. The weakened condition will necessarily make the proof longer and, therefore, is not considered in this paper.

Before we state the main results of this section, we need to introduce some notation that will be used throughout this paper. For sequences of real positive numbers $a_n$ and a sequence of real or random numbers $b_n$ that may depend on model parameters, $b_n = O_u(a_n)$ if, for any $n$, $P(|b_n| \leq Ma_n) = 1$, for some $0 < M < \infty$, which does not depend on parameters but could be random. That is, $b_n/a_n$ is bounded uniformly in the parameters. Similarly, we write $b_n = o_u(a_n)$ to mean that, with a probability 1, $b_n/a_n$ converges to 0 uniformly in parameters. The following theorem compares the tapered log-likelihood function with the untapered one, and their derivatives. This theorem is essential to the establishment of the asymptotic properties of the tapered MLE to be given in the subsequent theorem.

THEOREM 1. *Under the assumptions* (A1) *and* (A2), *uniformly in* $(\theta, \sigma^2) \in J$ *and with* $P_0$-*probability 1,*

$$l_{n,\text{tap}}(\theta, \sigma^2) = l_n(\theta, \sigma^2) + o_u(n^{1/2}), \quad (2.1)$$

$$\frac{\partial}{\partial \theta} l_{n,\text{tap}}(\theta, \sigma^2) = \frac{\partial}{\partial \theta} l_n(\theta, \sigma^2) + o_u(n^{1/2}), \quad (2.2)$$

*where* $P_0$ *is the probability measure corresponding to the true parameter values* $\sigma_0^2, \theta_0$.

The next theorem establishes the strong consistency and the asymptotic distribution of the tapered MLE. Comparing the asymptotic distribution of MLE of $\sigma_0^2 \theta_0$ in Ying (1993) and that in the following theorem, we see that the tapered MLE is asymptotically equally efficient.

THEOREM 2. *Assume* (A1) *and* (A2) *hold, and let* $(\hat{\theta}_{n,\text{tap}}, \hat{\sigma}_{n,\text{tap}}^2)$ *maximize the tapered likelihood function over* $(\theta, \sigma^2) \in J$. *Then, as* $n \to \infty$,

$$P_0\Big(\lim_{n \to \infty} \hat{\theta}_{n,\text{tap}} \hat{\sigma}_{n,\text{tap}}^2 = \theta_0 \sigma_0^2\Big) = 1, \tag{2.3}$$

$$\sqrt{n}(\hat{\theta}_{n,\text{tap}} \hat{\sigma}_{n,\text{tap}}^2 - \theta_0 \sigma_0^2) \xrightarrow{d} N(0, 2(\theta_0 \sigma_0^2)^2), \tag{2.4}$$

where $P_0$ is the probability measure corresponding to the true parameter values $\sigma_0^2, \theta_0$.

**3. General Matérn model.** In this section, we will focus on studying the asymptotics of tapered MLE for a general Matérn model. We assume the underlying process is stationary with mean 0 and the following isotropic Matérn covariogram:

$$K(h; \sigma^2, \theta, \nu) = \frac{\sigma^2 (\theta h)^\nu}{\Gamma(\nu) 2^{\nu-1}} \mathcal{K}_\nu(\theta h), \qquad h > 0, \tag{3.1}$$

with unknown $\sigma^2, \theta$ and known $\nu$, where $\mathcal{K}_\nu$ is the modified Bessel function of order $\nu$ [see Abramowitz and Stegun (1967), pages 375 and 376], $\sigma^2$ is the covariance parameter, $\theta$ is the scale parameter and $\nu$ is the smoothness parameter. Further, assume that the process is observed at $n$ sites $t_1, t_2, \ldots, t_n$ in a bounded interval $D \subset \mathbb{R}$, and write $\mathbf{X}_n = (X(t_1), \ldots, X(t_n))'$. Zhang (2004) noted that neither $\sigma^2$ or $\theta$ is consistently estimable under the fixed-domain asymptotic framework, but the quantity $\sigma^2 \theta^{2\nu}$ is consistently estimable. Furthermore, this consistently estimable quantity is more important to prediction than the parameters $\sigma^2$ and $\theta$.

The primary focus of this section is to establish the asymptotic distribution of the estimators for $\sigma^2 \theta^{2\nu}$. This is a more difficult problem than in the exponential case, and we cope with it by considering an easy version of the problem. Following Zhang (2004), we fix $\theta$ at an arbitrarily chosen value $\theta_1$ and consider the following estimators:

$$\hat{\sigma}_n^2 = \operatorname{Arg\,Max} l_n(\theta_1, \sigma^2), \tag{3.2}$$

$$\hat{\sigma}_{n,\text{tap}}^2 = \operatorname{Arg\,Max} l_{n,\text{tap}}(\theta_1, \sigma^2), \tag{3.3}$$

where $l_n(\theta_1, \sigma^2)$ and $l_{n,\text{tap}}(\theta_1, \sigma^2)$ are the log-likelihood function and the tapered log-likelihood function, respectively.

We make the following assumption on the spectral density of the taper $K_{\text{tap}}(h)$. Similar conditions were used in Furrer, Genton and Nychka (2006) and Kaufman, Schervish and Nychka (2008). Our condition here is stronger, and it is necessary for our approach to deriving the asymptotic distribution of tapered MLE:

(A3) The spectral density of the taper, denoted by $f_{\text{tap}}(\lambda)$, satisfies for some constant $\varepsilon > \max\{1/2, 1-\nu\}$ and $0 < M < \infty$

$$(3.4) \qquad f_{\text{tap}}(\lambda) \leq \frac{M}{(1+\lambda^2)^{\nu+1/2+\varepsilon}}.$$

We note that taper condition (3.4) is satisfied by some well-known tapers. For example, Wendland tapers (1995, 1998) have isotropic spectral densities that are continuous and satisfy $g_{d,k}(\lambda) \leq M(1+\lambda^2)^{-d/2-k-1/2}$ for some constant $M$, where $d$ is the dimension of the domain ($d=1$ in this work). Therefore, condition (3.4) is satisfied if $k > \max\{1/2, \nu\}$. Furrer, Genton and Nychka (2006) gave explicit tail limits for two Wendland tapers $K_1(h;\gamma) = (1-\frac{h}{\gamma})^4_+(1+4\frac{h}{\gamma}), \gamma > 0$ and $K_2(h;\gamma) = (1-\frac{h}{\gamma})^6_+(1+6\frac{h}{\gamma} + \frac{35h^2}{3\gamma^2}), \gamma > 0$ ($x_+ = \max\{0,x\}$), and showed that $\lambda^4 g_1(\lambda) \to 120/(\pi\gamma^3)$ and $\lambda^6 g_2(\lambda) \to 17920/(\pi\gamma^5)$, as $\lambda \to \infty$, where $g_i$ is the spectral density of $K_i$ ($i = 1, 2$). Therefore, condition (3.4) holds if $\nu < 1$ for taper $K_1$ and $\nu < 2$ for taper $K_2$.

One important probabilistic tool we will extensively use is the equivalence of probability measures. The assumption (A3) implies that the tapered covariance function specifies a Gaussian measure that is equivalent to the Gaussian measure specified by the true covariance function [Kaufman, Schervish and Nychka (2008)]. It readily follows that $\hat{\sigma}^2_{n,\text{tap}} \theta_1^{2\nu}$ is a strongly consistent estimator of $\sigma_0^2 \theta_0^{2\nu}$ [e.g., Kaufman, Schervish and Nychka (2008)].

The main results in this section are the following three theorems. The next theorem is a general result about two equivalent Gaussian measures and is not restricted to the case of covariance tapering. It will be used to prove the other two theorems.

THEOREM 3. *Let $X(t), t \in \mathbb{R}$ be a stationary Gaussian process having mean zero and an isotropic covariogram $K_j$ and a continuous spectral density $f_j$ under measure $P_j, j = 0, 1$. Assume the process is observed at $t_1, t_2, \ldots$ in a bounded interval $D$, and let $\mathbf{X}_n = (X(t_1), \ldots, X(t_n))'$. If*

$$(3.5) \qquad \liminf_{\lambda \to \infty} f_0(\lambda)|\lambda|^{r_1} > 0 \quad \text{and} \quad \limsup_{\lambda \to \infty} f_0(\lambda)|\lambda|^{r_1} < \infty$$

*and*

$$(3.6) \quad h(\lambda) = \frac{f_1(\lambda)}{f_0(\lambda)} - 1 = O(|\lambda|^{-r_2}), \qquad \lambda \to \infty \qquad \text{for some } r_2 > 1,$$

*then*

$$(3.7) \qquad E_0(\mathbf{X}'_n(\mathbf{V}^{-1}_{1,n} - \mathbf{V}^{-1}_{0,n})\mathbf{X}_n) = O(1),$$

*where $V_{j,n}$ is the covariance matrix of $\mathbf{X}_n$ given by the covariogram $K_j$, $j = 0, 1$, and $E_0$ is the expectation with respect to $P_0$.*

We note that condition (3.6) is stronger than the equivalence of the two Gaussian measures corresponding to the two spectral densities $f_0$ and $f_1$. Indeed, under condition (3.5), the two Gaussian measures are equivalent if (3.6) holds for some $r_2 > 1/2$. However, equivalence alone cannot imply (3.7), and we need stronger conditions than the equivalence of two measures. We will show later that condition (A3) implies that (3.6) holds, for some $r_2 > 1$, if $f_0$ and $f_1$ represent the spectral densities of the true and tapered covariograms, respectively.

THEOREM 4. *Suppose condition* (A3) *is satisfied, and the underlying process is stationary Gaussian having a mean 0 and a Matérn covariance function, and the sampling locations $\{t_1, t_2, \ldots\}$ are from a bounded interval. Then, for any fixed $\theta_1 > 0$, with $P_0$-probability 1, uniformly in $\sigma^2 \in [w, v]$,*

$$(3.8) \qquad l_{n,\text{tap}}(\theta_1, \sigma^2) = l_n(\theta_1, \sigma^2) + O_u(1),$$

$$(3.9) \qquad \frac{\partial}{\partial \sigma^2} l_{n,\text{tap}}(\theta_1, \sigma^2) = \frac{\partial}{\partial \sigma^2} l_n(\theta_1, \sigma^2) + O_u(1),$$

*where $P_0$ is the probability measure corresponding to the true parameter values $\sigma_0^2, \theta_0, \nu$.*

Next, we give the asymptotic distributions for both exact MLE and tapered MLE of the consistently estimable quantity $\sigma^2 \theta^{2\nu}$.

THEOREM 5. *Assume that the underlying process is stationary Gaussian having a mean 0 and a Matérn covariance function with a known smoothness parameter $\nu$, and the sampling locations $\{t_1, t_2, \ldots\}$ are from a bounded interval:*

(i) *For any fixed $\theta_1$,*

$$(3.10) \qquad \sqrt{n}(\hat{\sigma}_n^2 \theta_1^{2\nu} - \sigma_0^2 \theta_0^{2\nu}) \xrightarrow{d} N(0, 2(\sigma_0^2 \theta_0^{2\nu})^2).$$

(ii) *In addition, if the taper satisfies condition* (A3),

$$(3.11) \qquad \sqrt{n}(\hat{\sigma}_{n,\text{tap}}^2 \theta_1^{2\nu} - \sigma_0^2 \theta_0^{2\nu}) \xrightarrow{d} N(0, 2(\sigma_0^2 \theta_0^{2\nu})^2).$$

Theorem 5 implies that the covariance tapering does not reduce the asymptotic efficiency for the Matérn model. In this paper, we are not able to show that (3.11) remains true if $\theta_1$ is replaced by the MLE of $\theta$. Therefore, the results in this theorem are not as strong as those in Theorem 2. More work will need to be done to extend Theorem 2 to the general Matérn case.

We note that asymptotic distributional results about the microergodic parameter $\sigma^2 \theta^{2\nu}$ in the general Matérn class have not appeared in literature. Theorem 5(i) is the first of such results, and its proof requires a novel approach.

**4. Discussion.** There are some open problems for future research. First, for the Matérn model, the estimator of $\sigma^2\theta^{2\nu}$ is constructed by fixing $\theta$ at an arbitrary value. For a finite sample, common practice is to also estimate $\theta$. It is an interesting question to see if Theorem 5 still holds for the MLE $\hat{\sigma}^2\hat{\theta}^{2\nu}$ and the tapered MLE $\hat{\sigma}^2_{n,\text{tap}}\hat{\theta}^{2\nu}_{n,\text{tap}}$. Our conjecture is that Theorem 5 can be extended to this case.

The main results in Sections 2 and 3 are for the processes with one-dimensional index. It is a more interesting problem to study the high-dimensional case. However, our techniques in Section 3 cannot be extended to obtain analogous asymptotic distribution in the high-dimensional case. For example, for a $d$-dimensional process, we would need (3.6) to hold for some $r_2 > d$ in order for the proof to carry through. Unfortunately, for the Matérn model, (3.6) cannot hold for any $r_2 > 2$. The high-dimensional case calls for new techniques for establishing asymptotic distributions. A referee suggested letting the bandwidth $\gamma$ vary and go to 0 as $n$ increases to $\infty$. This is a natural scheme in the fixed-domain asymptotic framework. We believe that everything in Section 2 carries through if the bandwidth goes to 0 not too fast. When the bandwidth of the taper depends on $n$, it is not obvious if our techniques in Section 3 still apply, because the properties of equivalence of probability measures are no longer directly applicable.

## APPENDIX A: PROOFS FOR SECTION 2

In the sequel, we often suppress $n$ in the subscripts. For example, write $t_k = t_{k,n}$, $\Delta_k = \Delta_{k,n}$. We will need three lemmas for the proofs of the theorems in Section 2.

LEMMA 1. *Let $X(t)$ be the Gaussian Ornstein–Uhlenbeck process, and assume* (A1) *holds. Denote $E(X(t_i)|X(t_j), j\neq i) = -\sum_{j\neq i} b_{ij,n}(\theta)X(t_j)$, $1 \leq i \leq n$, $\text{Var}(X(t_i)|X(t_j), j\neq i) = d_{i,n}(\theta,\sigma^2)$, which is written as $d_i$ for short. Then, for $1 < i < n$,*

(A.1)
$$b_{ii-1,n}(\theta) = -\frac{e^{-\theta\Delta_i}(1-e^{-2\theta\Delta_{i+1}})}{1-e^{-2\theta(\Delta_i+\Delta_{i+1})}},$$
$$b_{ii+1,n}(\theta) = -\frac{e^{-\theta\Delta_{i+1}}(1-e^{-2\theta\Delta_i})}{1-e^{-2\theta(\Delta_i+\Delta_{i+1})}},$$

(A.2)
$$b_{12,n}(\theta) = -e^{-\theta\Delta_2}, \qquad b_{nn-1,n}(\theta) = -e^{-\theta\Delta_n},$$
$$b_{ij,n}(\theta) = 0 \quad \text{for } |i-j| > 1.$$

*In addition, uniformly in $(\theta,\sigma^2) \in J$, $1 < i < n$, $1 \leq j \leq n$,*

$$d_1 = 2\sigma^2\theta\Delta_2 + O_u\!\left(\frac{1}{n^2}\right), \qquad d_n = 2\sigma^2\theta\Delta_n + O_u\!\left(\frac{1}{n^2}\right),$$

(A.3)
$$d_i = \frac{2\sigma^2 \theta \Delta_i \Delta_{i+1}}{\Delta_i + \Delta_{i+1}} + O_u\left(\frac{1}{n^2}\right),$$

(A.4)
$$b'_{ij,n}(\theta) = O_u\left(\frac{1}{n^2}\right), \qquad b'_{1j,n}(\theta) = O_u\left(\frac{1}{n}\right),$$
$$b'_{nj,n}(\theta) = O_u\left(\frac{1}{n}\right), \qquad \frac{\partial}{\partial \theta} d_j^{-1} = O_u(n).$$

PROOF. Note that $E(X(t_i)|X(t_j), j \neq i) = -\sum_{j \neq i} b_{ij,n}(\theta) X(t_j), 1 \leq i \leq n$ if and only if

$$\mathrm{Cov}\left(X(t_i) + \sum_{k \neq i} b_{ik,n}(\theta) X(t_k), X(t_j)\right) = 0 \qquad \text{for any } j \neq i, j = 1, \ldots, n.$$

We therefore prove (A.1) and (A.2) by verifying that

(A.5)
$$\mathrm{Cov}(X(t_i) + b_{i,i-1} X(t_{i-1}) + b_{i,i+1} X(t_{i+1}), X(t_j)) = 0$$
$$\text{for any } j \neq i,$$

where we let $b_{10} = b_{n,n+1} = 0$. For $i = 1$ or $n$, (A.5) readily follows the stationarity and the Markovian property of the Ornstein–Uhlenbeck process. For $1 < i < n$, (A.5) holds, because, if $j \geq i+1$, the left-hand side of (A.5) equals

$$\sigma^2 e^{-\theta(t_j - t_{i+1})} \left( e^{-\theta \Delta_{i+1}} - \frac{e^{-\theta \Delta_i}(1 - e^{-2\theta \Delta_{i+1}})}{1 - e^{-2\theta(\Delta_i + \Delta_{i+1})}} e^{-\theta(\Delta_i + \Delta_{i+1})} \right.$$
$$\left. - \frac{e^{-\theta \Delta_{i+1}}(1 - e^{-2\theta \Delta_i})}{1 - e^{-2\theta(\Delta_i + \Delta_{i+1})}} \right),$$

which is zero. We can get the similar expression when $j \leq i - 1$. Therefore, (A.5) is proved. Since $d_i = E(X(t_i) + b_{i,i-1} X(t_{i-1}) + b_{i,i+1} X(t_{i+1}))^2$, straightforward calculation yields

$$d_1 = \sigma^2(1 - e^{-2\theta \Delta_2}), \qquad d_n = \sigma^2(1 - e^{-2\theta \Delta_n}),$$
$$d_i = \sigma^2 \frac{(1 - e^{-2\theta \Delta_i})(1 - e^{-2\theta \Delta_{i+1}})}{1 - e^{-2\theta(\Delta_i + \Delta_{i+1})}}, \qquad 1 < i < n.$$

Then, (A.3) follows the Taylor expansion. To establish the properties of the derivatives in (A.4), we repeatedly use the Taylor expansion. Here we only provide a proof for $b'_{ij,n}(\theta) = O_u(1/n^2)$ for $1 < i < n$, since all other derivatives can be proved similarly. Since $b_{ij} = 0$ if $|i - j| > 1$, we only need to consider $j = i - 1$ or $i + 1$.

Write the derivative

$$b'_{ii-1,n}(\theta) = A/(1 - e^{-2\theta(\Delta_i + \Delta_{i+1})})^2,$$

where

$$A = (\Delta_i e^{-\theta\Delta_i} - (\Delta_i + 2\Delta_{i+1})e^{-\theta\Delta_i - 2\theta\Delta_{i+1}})(1 - e^{-2\theta(\Delta_i + \Delta_{i+1})})$$
$$- (-e^{-\theta\Delta_i} + e^{-\theta\Delta_i - 2\theta\Delta_{i+1}})2(\Delta_i + \Delta_{i+1})e^{-2\theta(\Delta_i + \Delta_{i+1})}$$
(A.6)
$$= \Delta_i e^{-\theta\Delta_i} - (\Delta_i + 2\Delta_{i+1})(e^{-\theta\Delta_i - 2\theta\Delta_{i+1}} - e^{-3\theta\Delta_i - 2\theta\Delta_{i+1}})$$
$$- \Delta_i e^{-3\theta\Delta_i - 4\theta\Delta_{i+1}}.$$

Note that

$$\left| \frac{1}{1 - e^{-2\theta(\Delta_i + \Delta_{i+1})}} - \frac{1}{2\theta(\Delta_i + \Delta_{i+1})} \right|$$

is uniformly bounded and $n(\Delta_i + \Delta_{i+1})$ is bounded away from 0 and $\infty$ by assumption (A1). Hence, $1/(1 - e^{-2\theta(\Delta_i + \Delta_{i+1})}) = O(1/n)$, and it suffices to show that $A$ is $O_u(1/n^4)$. Using, again, the fact that $\Delta_i = O_u(1/n)$ and applying the Taylor expansion, we get

$$\Delta_i e^{-\theta\Delta_i} = \Delta_i - \theta\Delta_i^2 + (1/2)\theta^2\Delta_i^3 + O_u(1/n^4),$$
$$- (\Delta_i + 2\Delta_{i+1})(e^{-\theta\Delta_i - 2\theta\Delta_{i+1}} - e^{-3\theta\Delta_i - 2\theta\Delta_{i+1}})$$
$$= -2\theta\Delta_i^2 + 4\theta^2\Delta_i^3 + 12\theta^2\Delta_i^2\Delta_{i+1}$$
$$- 4\theta\Delta_i\Delta_{i+1} + 8\theta^2\Delta_i\Delta_{i+1}^2 + O_u(1/n^4),$$
$$-\Delta_i e^{-3\theta\Delta_i - 4\theta\Delta_{i+1}} = -\Delta_i + 3\theta\Delta_i^2 + 4\theta\Delta_i\Delta_{i+1} - (9/2)\theta^2\Delta_i^3$$
$$- 12\theta^2\Delta_i^2\Delta_{i+1} - 8\theta^2\Delta_i\Delta_{i+1}^2 + O_u(1/n^4).$$

All the terms except $O_u(1/n^4)$ are canceled out. Therefore, $A = O_u(1/n^4)$ and $b'_{ii-1,n}(\theta) = O_u(1/n^2)$. Similarly, we can show $b'_{ii+1,n}(\theta) = O_u(1/n^2)$. □

We now introduce the following notations. Let $\tilde{\mathbf{O}}_n$ denote a matrix of which the elements are $O_u(1/n)$ except those in the first and last rows, which are uniformly bounded; that is, $O_u(1)$. Denote, by $\check{\mathbf{O}}_n$, the matrix whose $(i,j)$th element is $O_u(1)$ if $i = 1$ or $n$ or $i = j$, and is $O_u(1/n)$ otherwise. Therefore,

$$\tilde{\mathbf{O}}_n = \begin{pmatrix} O_u(1) & \cdots & O_u(1) \\ O_u(1/n) & \cdots & O_u(1/n) \\ \cdots & \cdots & \cdots \\ O_u(1/n) & \cdots & O_u(1/n) \\ O_u(1) & \cdots & O_u(1) \end{pmatrix},$$

$$\check{\mathbf{O}}_n = \tilde{\mathbf{O}}_n + \begin{pmatrix} O_u(1) & & & & \\ & O_u(1) & & & \\ & & \ddots & & \\ & & & O_u(1) & \\ & & & & O_u(1) \end{pmatrix}.$$

LEMMA 2. *Under assumptions* (A1) *and* (A2), *uniformly in* $\theta \in [a, b]$:

(i) $\mathbf{V}_n^{-1}(\mathbf{V}_n \circ \mathbf{T}_n) = \mathbf{I}_n + \tilde{\mathbf{O}}_n, \qquad \mathbf{V}_n^{-1}\dfrac{\partial \mathbf{V}_n}{\partial \theta} = \breve{\mathbf{O}}_n,$

(ii) $\dfrac{\partial}{\partial \theta}(\mathbf{V}_n^{-1}(\mathbf{V}_n \circ \mathbf{T}_n)) = \tilde{\mathbf{O}}_n, \qquad \dfrac{\partial}{\partial \theta}\left(\mathbf{V}_n^{-1}\dfrac{\partial \mathbf{V}_n}{\partial \theta}\right) = \breve{\mathbf{O}}_n,$

(iii) $1 < \det(\mathbf{V}_n^{-1}(\mathbf{V}_n \circ \mathbf{T}_n)) = O_u(1), \qquad (\mathbf{V}_n^{-1}(\mathbf{V}_n \circ \mathbf{T}_n))^{-1} = \mathbf{I}_n + \tilde{\mathbf{O}}_n,$

where $\mathbf{I}_n$ is the $n \times n$ identity matrix.

From the definitions of $\tilde{\mathbf{O}}_n$ and $\breve{\mathbf{O}}_n$, we have

(A.7) $\qquad \tilde{\mathbf{O}}_n \breve{\mathbf{O}}_n = \tilde{\mathbf{O}}_n, \qquad \breve{\mathbf{O}}_n \tilde{\mathbf{O}}_n = \tilde{\mathbf{O}}_n, \qquad \tilde{\mathbf{O}}_n \tilde{\mathbf{O}}_n = \tilde{\mathbf{O}}_n.$

Then, Lemma 2(i) and (ii) imply

(A.8) $\qquad \dfrac{\partial}{\partial \theta}(\mathbf{V}_n^{-1}(\mathbf{V}_n \circ \mathbf{T}_n))^{-1} = \tilde{\mathbf{O}}_n.$

PROOF OF LEMMA 2. We can assume $\sigma^2 = 1$ without loss of any generality, because all quantities in the lemma do not depend on $\sigma^2$. We will repeatedly use Lemma 1 and particularly the fact that $\mathbf{V}_n^{-1}$ is a band matrix. The proof involves tedious computation, and we will keep a balance between brevity and clarity.

Several quantities in the lemma are of the form $\mathbf{V}_n^{-1}(\mathbf{V}_n \circ \mathbf{Q})$, where $\mathbf{Q}$ is an $n \times n$ matrix whose $(i, j)$th element is $\varrho(t_i - t_j)$ for some even function $\varrho(t)$ that has a bounded second derivative on $[-1, 0) \cup (0, 1]$. If the limits of the derivative $\varrho'(0+) = \lim_{t \to +} \varrho'(t)$ and $\varrho'(0-) = \lim_{t \to 0-} \varrho'(t)$ exist and are finite, we show now

(A.9) $\mathbf{V}_n^{-1}(\mathbf{V}_n \circ \mathbf{Q}) = \varrho(0)\mathbf{I}_n + (\varrho'(0+) - \varrho'(0-))$
$$\times \operatorname{diag}\{O_u(1), \ldots, O_u(1)\} + \tilde{\mathbf{O}}_n,$$

where $\operatorname{diag}(O_u(1), \ldots, O_u(1))$ denotes a diagonal $n \times n$ matrix with bounded elements. There are immediate corollaries from (A.9). First, it implies $\mathbf{V}_n^{-1}(\mathbf{V}_n \circ \mathbf{Q}) = \breve{\mathbf{O}}_n$. Second, by taking $\varrho(t) = -|t|$, we get $\mathbf{V}_n^{-1}(\partial \mathbf{V}_n / \partial \theta) = \breve{\mathbf{O}}_n$. Last, if $\varrho(t) = K_{\mathrm{tap}}(|t|)$, then $\varrho'(0+) = \varrho'(0-)$ and $\mathbf{V}_n^{-1}(\mathbf{V}_n \circ \mathbf{T}_n) = \mathbf{I}_n + \tilde{\mathbf{O}}_n$ because $K_{\mathrm{tap}}(0) = 1$.

To prove (A.9), let $\omega_{ij}$ denote the $(i, j)$th element of $\mathbf{V}_n^{-1}(\mathbf{V}_n \circ \mathbf{Q})$. Let $b_{ij}$ be defined in Lemma 1. Hereafter, the parameter and subscript $n$ are suppressed. Then, it is well known that the $(i, j)$th element of $\mathbf{V}_n^{-1}$ is $b_{ij}/d_i$. Write $b_{ij} = 0$ if $j < 1$ or $j > n$ and $t_0 = t_1, t_{n+1} = t_n$. Since $b_{ij} = 0$ if $|i - j| > 1$,

(A.10) $\qquad \omega_{ij} = d_i^{-1} \sum_{k=i-1}^{i+1} b_{ik} K(|t_k - t_j|) \varrho(t_k - t_j).$

For any $i > j$, and $k = i - 1$ or $i + 1$, we have $t_k - t_j \geq 0$. Hence, the Taylor theorem implies

$$\varrho(t_k - t_j) = \varrho(t_i - t_j) + \varrho'(t_i - t_j)(t_k - t_i) + \varrho''(t_i - t_j + \xi(t_k - t_j))(t_k - t_i)^2/2,$$

for some $\xi \in (0, 1)$. Since $\varrho$ has a bounded second derivative on $(0, 1)$, and $t_k - t_i = O(1/n)$, we have

(A.11) $\quad \varrho(t_k - t_j) = \varrho(t_i - t_j) + \varrho'(t_i - t_j)(t_k - t_i) + O(1/n^2).$

Then,

(A.12)
$$\omega_{ij} = d_i^{-1} \varrho(t_i - t_j) \sum_{k=i-1}^{i+1} b_{ik} K(t_k - t_j)$$
$$+ d_i^{-1} \varrho'(t_i - t_j) \sum_{k=i-1}^{i+1} K(t_k - t_j)(t_k - t_i) b_{ik} + O_u(1/n),$$

where we have used $d_i^{-1} = O_u(n)$. Note that $d_i^{-1} \sum_{k=i-1}^{i+1} b_{ik} K(t_k - t_j)$ is the $(i, j)$ element of $\mathbf{D}_n^{-1} \mathbf{B}_n \mathbf{V}_n = \mathbf{I}_n$. Hence, the first summand in (A.12) equals $\varrho(0) 1_{\{i=j\}}$.

Similar to the establishment of (A.11), we can show

$$K(t_k - t_j) = K(t_i - t_j) + K'(t_i - t_j)(t_k - t_i) + O_u(1/n^2).$$

It follows that, for $i > j$,

(A.13) $\quad \omega_{ij} = d_i^{-1} \varrho'(t_i - t_j) K(t_i - t_j) \sum_{k=i-1}^{i+1} (t_k - t_i) b_{ik} + O_u(1/n).$

By utilizing the explicit expressions of $b_{ij}$ given in Lemma 1, we can show

(A.14)
$$\sum_{k=i-1}^{i+1} (t_k - t_i) b_{ik} = b_{i,i+1} \Delta_{i+1} - b_{i,i-1} \Delta_i$$
$$= \begin{cases} O_u(1/n^2), & \text{if } 1 < i < n, \\ O_u(1/n), & \text{if } i = 1 \text{ or } n. \end{cases}$$

Then, for $i > j$

(A.15) $\quad \omega_{ij} = \begin{cases} O_u(1), & \text{if } i = 1 \text{ or } n, \\ O_u(1/n), & \text{if } 1 < i < n. \end{cases}$

In view of the fact that $\varrho$ is an even function, we can show, similarly, that (A.15) holds for $i < j$.

Now, let us consider $w_{ii}$. First, note that

(A.16) $\quad \varrho(t_{i-1} - t_i) = \varrho(0) + \varrho'(0-)(t_{i-1} - t_i) + O(1/n^2),$

(A.17) $\quad \varrho(t_{i+1} - t_i) = \varrho(0) + \varrho'(0+)(t_{i+1} - t_i) + O(1/n^2).$

Since $d_i^{-1} \sum_{k=i-1}^{i+1} b_{ik} K(t_k - t_i) = 1$,

$$\omega_{ii} = d_i^{-1} \sum_{k=i-1}^{i+1} b_{ik} K(t_k - t_i) \varrho(t_k - t_i)$$

$$= \varrho(0) + d_i^{-1} \{ b_{i,i-1} K(t_{i-1} - t_i) \varrho'(0-)(t_{i-1} - t_i)$$
$$+ b_{i,i+1} K(t_{i+1} - t_i) \varrho'(0+)(t_{i+1} - t_i) \} + O_u(1/n^2).$$

Since $K(h) = K(0) + K'(0)h + o_u(h)$ as $h \to 0$,

$$\omega_{ii} = \varrho(0) + K(0) d_i^{-1} \{ b_{i,i-1} \varrho'(0-)(t_{i-1} - t_i) + b_{i,i+1} \varrho'(0+)(t_{i+1} - t_i) \}$$
$$+ O_u(1/n^2),$$

which can be rewritten as

$$\omega_{ii} = \varrho(0) + \varrho'(0-)K(0) d_i^{-1} \sum_{k=i-1}^{i+1} (t_k - t_i) b_{ik}$$

(A.18)
$$+ [\varrho'(0+) - \varrho'(0-)] K(0) d_i^{-1} b_{i,i+1} \Delta_{i+1} + O_u(1/n^2).$$

Then, (A.9) follows from (A.14), (A.15) and (A.18). (i) is therefore proved.

To prove (ii), we will use the following well-known fact [see, e.g., Graybill (1983), pages 357 and 358]:

(A.19) $$\frac{\partial}{\partial \theta} \mathbf{V}_n^{-1} = -\mathbf{V}_n^{-1} \frac{\partial \mathbf{V}_n}{\partial \theta} \mathbf{V}_n^{-1}.$$

Applying Lemma 2(i),

$$\frac{\partial}{\partial \theta} (\mathbf{V}_n^{-1}(\mathbf{V}_n \circ \mathbf{T}_n)) = -\mathbf{V}_n^{-1} \frac{\partial \mathbf{V}_n}{\partial \theta} \mathbf{V}_n^{-1}(\mathbf{V}_n \circ \mathbf{T}_n) + \mathbf{V}_n^{-1} \left( \frac{\partial \mathbf{V}_n}{\partial \theta} \circ \mathbf{T}_n \right)$$

$$= -\mathbf{V}_n^{-1} \frac{\partial \mathbf{V}_n}{\partial \theta} (\mathbf{I}_n + \tilde{\mathbf{O}}_n) + \mathbf{V}_n^{-1} \left( \frac{\partial \mathbf{V}_n}{\partial \theta} \circ \mathbf{T}_n \right)$$

$$= \mathbf{V}_n^{-1} \left( \frac{\partial \mathbf{V}_n}{\partial \theta} \circ (\mathbf{T}_n - \mathbf{J}_n) \right) + \tilde{\mathbf{O}}_n,$$

which is clearly $\tilde{\mathbf{O}}_n$ from (A.9) by taking $\varrho(t) = -|t|(K_{\text{tap}}(|t|) - 1)$ that is differentiable at 0, where $\mathbf{J}_n$ is a matrix of all 1's.

Next, we will show $\frac{\partial}{\partial \theta}(\mathbf{V}_n^{-1} \frac{\partial \mathbf{V}_n}{\partial \theta}) = \check{\mathbf{O}}_n$ similarly. Write

(A.20) $$\frac{\partial}{\partial \theta}\left( \mathbf{V}_n^{-1} \frac{\partial \mathbf{V}_n}{\partial \theta} \right) = -\mathbf{V}_n^{-1} \frac{\partial \mathbf{V}_n}{\partial \theta} \mathbf{V}_n^{-1} \frac{\partial \mathbf{V}_n}{\partial \theta} + \mathbf{V}_n^{-1} \frac{\partial^2 \mathbf{V}_n}{\partial \theta^2}.$$

By (i), the first term on the right-hand side of (A.20) is $\check{\mathbf{O}}_n \check{\mathbf{O}}_n = \check{\mathbf{O}}_n$, and the second term is $\tilde{\mathbf{O}}_n$, because $\mathbf{V}_n^{-1} \frac{\partial^2 \mathbf{V}_n}{\partial \theta^2} = \mathbf{V}_n^{-1}(\mathbf{V}_n \circ \mathbf{Q})$ with $\varrho(t) = t^2$, which has a continuous second derivative so that (A.9) applies. This completes the proof of Lemma 2(ii).

Let $\mathbf{A}_n = \mathbf{V}_n^{-1}(\mathbf{V}_n \circ \mathbf{T}_n)$ and $a_{ij}$ denote the $(i,j)$th element of $\mathbf{A}_n$. We now apply a series of column operations, so that $\mathbf{A}_n$ becomes $\mathbf{I}_n + \mathbf{\Omega}_n$ and each of the operations retains the determinant of $\mathbf{A}_n$, where $\mathbf{\Omega}_n$ is a matrix whose elements are bounded by $M/n$ for some constant $M$ not depending on $\theta$; that is, $\mathbf{\Omega}_n(i,j) \leq M/n$. We have shown that $\mathbf{A}_n = \mathbf{I}_n + \tilde{\mathbf{O}}_n$, where elements of $\tilde{\mathbf{O}}_n$ are $O_u(1/n)$ except those that are on the first and last rows that are bounded. We can subtract from the $j$th column the first column multiplied by the $(1,j)$th element of $\mathbf{A}_n$, $2 < j < n$. Then, all elements in the first row are $O_u(1/n)$, except the $(1,1)$th element, which is $1 + O_u(1/n)$ and remains unchanged throughout the operations. Similarly, we can reduce the elements in the last row to $O_u(1/n)$ except the last $(n,n)$th element. Applying the Hadamard inequality [Bellman (1970), page 130], we can show there exists some constant $M$ such that

$$\det(\mathbf{A}_n) = \det(\mathbf{I}_n + \mathbf{\Omega}_n) \leq ((1 + M/n)^2 + (n-1)(M/n)^2)^{n/2},$$

which is bounded.

To show $\mathbf{A}_n^{-1} = \mathbf{I}_n + \tilde{\mathbf{O}}_n$, we first note that by Oppenheim's inequality [Mirsky (1955), page 421], which yields the inequality for the determinant of Hadamard product of positive definite matrices, $\det(\mathbf{V}_n \circ \mathbf{T}_n) > \det(\mathbf{V}_n) \prod_{1 \leq i \leq n} t_{ii}$ where $t_{ii}$ is the diagonal element of $\mathbf{T}_n$. Therefore, $\det(\mathbf{A}_n) > 1$. We only need to show that the $(i,j)$th cofactor

$$A_{ij} = \det(\mathbf{A}_n) 1_{\{i=j\}} + O_u(1/n) + 1_{\{j=1 \text{ or } j=n\}} O_u(1).$$

Similar to proving $\det(\mathbf{A}_n) = O_u(1)$, we can show all the $(n-1)$ by $(n-1)$ cofactors are also $O_u(1)$. In addition, $A_{ij} = O_u(1/n)$ for $1 < j < n, i \neq j$ since it has one row of elements $O_u(1/n)$ and replacing that row with $O_u(1)$ would yield a bounded determinant. To complete proof of the lemma, it remains to show $A_{ii} = \det(\mathbf{A}_n) + O_u(1/n), 1 < i < n$, which is true by Laplace expansion $\det(\mathbf{A}_n) = (1 + O_u(1/n))A_{ii} + \sum_{j : j \neq i} a_{ij} A_{ij}$ and observing that $\sum_{j : j \neq i} a_{ij} A_{ij} = O_u(1/n)$, for $1 < i < n$. □

LEMMA 3. *For any $\theta \in [a,b]$, let $S_n(\theta)$, $n = 1, 2, \ldots$, be a sequence of random variables such that $E(S_n(\theta)) = O_u((\log n)^r)$, $E[S_n(\theta) - ES_n(\theta)]^6 = O_u((\log n)^r)$ uniformly in $\theta$ for some constant $r > 0$. Assume that, with probability one, $S_n(\theta)$ is differentiable with respect to $\theta$ and $S'_n(\theta) = O_u(n^2(\log n)^r)$ uniformly in $\theta$. Then,*

$$\sup_{\theta \in [a,b]} |S_n(\theta)| = o(n^{1/2}) \quad a.s.$$

For the ease of notation, we will suppress the dependence of any quantity on $n$ and parameters [e.g., $b_{ij,n}(\theta) = b_{ij}$, $d_{i,n}(\theta, \sigma^2) = d_i$], wherever confusion does not arise, throughout the rest of the paper.

PROOF OF LEMMA 3. Let $a = \theta_0 < \theta_1 < \cdots < \theta_{M_n} = b$ partition $[a, b]$ into intervals of equal length, where $M_n$ is the integer part of $n^{3/2+\alpha}$ for some $0 < \alpha < 1/14$. Then,

$$(A.21) \quad \sup_{\theta \in [a,b]} |S_n(\theta)| \leq \max_{1 \leq k \leq M_n} |S_n(\theta_k)| + \max_{1 \leq k \leq M_n} \sup_{\theta \in [\theta_{k-1}, \theta_k]} |S_n(\theta_k) - S_n(\theta)|.$$

Because there exists constant $C > 0$ such that the sixth central moment of $S_n(\theta)$ is uniformly bounded by $C(\log n)^r$,

$$P\Big(\max_{1 \leq k \leq M_n} |S_n(\theta_k) - E(S_n(\theta_k))| \geq n^{1/2-\alpha}\Big)$$
$$\leq \sum_{k=1}^{M_n} P(|S_n(\theta_k) - E(S_n(\theta_k))| \geq n^{1/2-\alpha})$$
$$\leq M_n \frac{C(\log n)^r}{n^{3-6\alpha}} = \frac{C(\log n)^r}{n^{3/2-7\alpha}}.$$

Since $3/2 - 7\alpha > 1$ with $\alpha < 1/14$, it follows from Borel–Cantelli lemma that

$$\max_{1 \leq k \leq M_n} |S_n(\theta_k) - E(S_n(\theta_k))| = O(n^{1/2-\alpha}) \quad \text{a.s.}$$

Since $E(S_n(\theta)) = O((\log n)^r)$ uniformly in $\theta$, and

$$\max_{1 \leq k \leq M_n} |S_n(\theta_k)| \leq \max_{1 \leq k \leq M_n} |S_n(\theta_k) - E(S_n(\theta_k))| + \max_{1 \leq k \leq M_n} |E(S_n(\theta_k))|,$$

then

$$(A.22) \quad \max_{1 \leq k \leq M_n} |S_n(\theta_k)| = O(n^{1/2-\alpha}) \quad \text{a.s.}$$

On the other hand, for $\theta \in [\theta_{k-1}, \theta_k]$,

$$|S_n(\theta_k) - S_n(\theta)| \leq \sup_{\theta \in [a,b]} |S'_n(\theta)|(\theta_k - \theta_{k-1}) = O_u(n^{1/2-\alpha}(\log n)^r) \quad \text{a.s.}$$

Therefore,

$$(A.23) \quad \max_{1 \leq k \leq M_n} \sup_{\theta \in [\theta_{k-1}, \theta_k]} |S_n(\theta_k) - S_n(\theta)| = o(n^{1/2}) \quad \text{a.s.}$$

The proof is completed by combining (A.21), (A.22) and (A.23). □

PROOF OF THEOREM 1. Recall that the tapered and untapered log-likelihoods are given by (1.3) and (1.1), respectively. The proof of (2.1) consists of direct comparisons of the log determinants and the two quadratic forms. First, Lemma 2(iii) implies that

$$(A.24) \quad \log[\det(\mathbf{V}_n \circ \mathbf{T}_n)] = \log[\det(\mathbf{V}_n)] + O_u(1).$$

Define $\mathbf{H}_n(\theta) = (\mathbf{V}_n^{-1}(\mathbf{V}_n \circ \mathbf{T}_n))^{-1} - \mathbf{I}_n$. Then, $\mathbf{H}_n = \tilde{\mathbf{O}}_n$ by Lemma 2(iii). Because

$$(\mathbf{V}_n \circ \mathbf{T}_n)^{-1} = \mathbf{V}_n^{-1} + \mathbf{H}_n \mathbf{V}_n^{-1}, \tag{A.25}$$

$$\mathbf{X}_n'(\mathbf{V}_n \circ \mathbf{T}_n)^{-1}\mathbf{X}_n = \mathbf{X}_n'\mathbf{V}_n^{-1}\mathbf{X}_n + \mathbf{X}_n'\mathbf{H}_n \mathbf{V}_n^{-1}\mathbf{X}_n. \tag{A.26}$$

Proof of (2.1) would be completed if, uniformly in $(\theta, \sigma^2)$, with probability 1,

$$\mathbf{X}_n'\mathbf{H}_n \mathbf{V}_n^{-1}\mathbf{X}_n = o_u(n^{1/2}). \tag{A.27}$$

We will apply Lemma 3 to prove (A.27). Define

$$S_n(\theta) = \sigma^2 \mathbf{X}_n'\mathbf{H}_n \mathbf{V}_n^{-1}\mathbf{X}_n,$$

and note that $S_n(\theta)$ depends on $\theta$ but not on $\sigma^2$. In view of symmetry of $\mathbf{H}_n\mathbf{V}_n^{-1}$ by (A.25), we can write

$$E_0 S_n(\theta) = \sigma^2 \operatorname{trace}\{\mathbf{H}_n \mathbf{V}_n^{-1} \mathbf{V}_{n,0}\}, \tag{A.28}$$

where hereafter in the proof the expectation is evaluated under the true parameter $\sigma_0^2$ and $\theta_0$, and $\mathbf{V}_{n,0} = \mathbf{V}_n(\theta_0, \sigma_0^2)$. The $r$th cumulant of $\mathbf{X}_n'\mathbf{H}_n\mathbf{V}_n^{-1}\mathbf{X}_n$

$$\kappa_r = 2^{r-1}(r-1)! \operatorname{trace}\{\mathbf{H}_n \mathbf{V}_n^{-1} \mathbf{V}_{n,0}\}^r, \tag{A.29}$$

$r = 1, 2, \ldots$ [see Searle (1971), Theorem 1, page 55].

Next, we show that

$$\mathbf{V}_n^{-1}(\theta, \sigma^2)\mathbf{V}_n(\theta_0, \sigma_0^2) = O_u(1)\mathbf{I}_n + \tilde{\mathbf{O}}_n. \tag{A.30}$$

Then, it follows from (A.28)–(A.30) and (A.7) that the first moment and the sixth central moment of $S_n(\theta)$ are uniformly bounded, because the sixth central moment of $S_n(\theta)$ is $\kappa_6 + 15\kappa_4\kappa_2 + 10\kappa_3^2 + 15\kappa_2^2$, which is uniformly bounded because all of the four cumulants involved are uniformly bounded.

We now give explicit expression for the elements of $\mathbf{V}_n^{-1}(\theta, \sigma^2)$ based on Lemma 1 and the following well-known result [e.g., Ripley (1981), page 89]:

$$\mathbf{V}_n^{-1}(\theta, \sigma^2) = \mathbf{D}_n^{-1}(\theta, \sigma^2)\mathbf{B}_n(\theta), \tag{A.31}$$

where $\mathbf{B}_n(\theta) = (b_{ij,n}(\theta))_{1 \leq i,j \leq n}$ and $\mathbf{D}_n(\theta, \sigma^2) = \operatorname{diag}\{d_i(\theta, \sigma^2), i = 1, \ldots, n\}$, in which $b_{ij,n}(\theta)$, $d_i(\theta, \sigma^2)$ are defined as in Lemma 1 and $b_{ii,n}(\theta) = 1$.

For brevity, we drop the parameters in the matrices and write $\mathbf{B}_{n,0} = \mathbf{B}_n(\theta_0)$ and $\mathbf{D}_{n,0} = \mathbf{D}_n(\theta_0, \sigma_0^2)$. Decompose $\mathbf{V}_n^{-1}$ into

$$\mathbf{V}_n^{-1} = \mathbf{D}_n^{-1}\mathbf{B}_{n,0} + \mathbf{D}_n^{-1}(\mathbf{B}_{n,0} - \mathbf{B}_n) = \mathbf{A}_1 + \mathbf{A}_2.$$

Then,

$$\mathbf{A}_1 \mathbf{V}_{n,0} = \mathbf{D}_n^{-1}\mathbf{B}_{n,0}\mathbf{B}_{n,0}^{-1}\mathbf{D}_{n,0} = \mathbf{D}_n^{-1}\mathbf{D}_{n,0}$$

and the diagonals of $\mathbf{D}_n^{-1}\mathbf{D}_{n,0}$ converge uniformly to $(\sigma_0^2\theta_0)/(\sigma^2\theta)$ by (A.3). Therefore,

$$\text{(A.32)} \qquad \mathbf{A}_1\mathbf{V}_{n,0} = \text{diag}(O_u(1),\ldots,O_u(1)).$$

In addition,

$$\text{(A.33)} \qquad \mathbf{A}_2\mathbf{V}_{n,0} = \tilde{\mathbf{O}}_n \qquad \text{uniformly in } \theta \in [a,b].$$

Indeed, the absolute value of the $(i,j)$th element of $\mathbf{A}_2\mathbf{V}_{n,0}$, $1 < i < n$ is

$$\left|\sum_{k=1}^{n} d_i^{-1}(b_{ik,n}(\theta) - b_{ik,n}(\theta_0))\sigma_0^2 e^{-\theta_0|t_k-t_j|}\right|$$

$$\leq d_i^{-1}\sigma_0^2 \sum_{|k-i|\leq 1} |b_{ik,n}(\theta) - b_{ik,n}(\theta_0)|$$

$$= O_u\left(\frac{1}{n}\right),$$

where the last equality follows from (A.3), (A.4) and the Taylor theorem. Similarly, we can show the elements on the first and last rows are $O_u(1)$. Hence, (A.30) follows from (A.32), (A.33) immediately. Last, note that $\frac{\partial}{\partial\theta}S_n(\theta) = O_u(n^2)$ by Lemma 2. The conditions of Lemma 3 are satisfied. Therefore,

$$\text{(A.34)} \qquad \sup_{\theta \in [a,b]} S_n(\theta) = o(n^{1/2}),$$

which implies (A.27).

We have now proved (2.1). (2.2) can be proved similarly, and the remaining proof will be brief. The derivatives of the log likelihood functions can be written as

$$\text{(A.35)} \qquad \frac{\partial}{\partial\theta}l_n(\theta,\sigma^2) = -\text{trace}\left\{\mathbf{V}_n^{-1}\frac{\partial\mathbf{V}_n}{\partial\theta}\right\} + \mathbf{X}_n'\mathbf{V}_n^{-1}\frac{\partial\mathbf{V}_n}{\partial\theta}\mathbf{V}_n^{-1}\mathbf{X}_n,$$

$$\text{(A.36)} \qquad \frac{\partial}{\partial\theta}l_{n,\text{tap}}(\theta,\sigma^2) = -\text{trace}\left\{(\mathbf{V}_n\circ\mathbf{T}_n)^{-1}\left(\frac{\partial\mathbf{V}_n}{\partial\theta}\circ\mathbf{T}_n\right)\right\}$$
$$+ \mathbf{X}_n'(\mathbf{V}_n\circ\mathbf{T}_n)^{-1}\left(\frac{\partial\mathbf{V}_n}{\partial\theta}\circ\mathbf{T}_n\right)(\mathbf{V}_n\circ\mathbf{T}_n)^{-1}\mathbf{X}_n.$$

We first show that the two traces differ by $O_u(1)$. Write $\mathbf{A}_n = \mathbf{V}_n^{-1}(\mathbf{V}_n\circ\mathbf{T}_n)$. It is straightforward to verify

$$(\mathbf{V}_n\circ\mathbf{T}_n)^{-1}\left(\frac{\partial\mathbf{V}_n}{\partial\theta}\circ\mathbf{T}_n\right) = \mathbf{A}_n^{-1}\mathbf{V}_n^{-1}\frac{\partial\mathbf{V}_n}{\partial\theta}\mathbf{A}_n + \mathbf{A}_n^{-1}\frac{\partial\mathbf{A}_n}{\partial\theta}.$$

Then,

$$\text{trace}\left\{(\mathbf{V}_n\circ\mathbf{T}_n)^{-1}\left(\frac{\partial\mathbf{V}_n}{\partial\theta}\circ\mathbf{T}_n\right)\right\} = \text{trace}\left\{\mathbf{V}_n^{-1}\frac{\partial\mathbf{V}_n}{\partial\theta}\right\} + \text{trace}\left(\mathbf{A}_n^{-1}\frac{\partial\mathbf{A}_n}{\partial\theta}\right),$$

where the second trace in the right-hand side is clearly uniformly bounded by Lemma 2. Similarly, we can write

$$(\mathbf{V}_n \circ \mathbf{T}_n)^{-1}\left(\frac{\partial \mathbf{V}_n}{\partial \theta} \circ \mathbf{T}_n\right)(\mathbf{V}_n \circ \mathbf{T}_n)^{-1}$$

$$= \mathbf{V}_n^{-1}\frac{\partial \mathbf{V}_n}{\partial \theta}\mathbf{V}_n^{-1} + \mathbf{W}_n \mathbf{V}_n^{-1}$$

for some matrix $\mathbf{W}_n$, which is $\tilde{\mathbf{O}}_n$. Using the exact same technique for deriving (A.27), we can show

$$\mathbf{X}'_n \mathbf{W}_n \mathbf{V}_n^{-1} \mathbf{X}_n = o_u(n^{1/2}).$$

The proof is complete. □

PROOF OF THEOREM 2. First, for (2.3), it suffices to show that, for any $\varepsilon > 0$,

$$(A.37) \quad P_0\Big(\inf_{\{(\theta,\sigma^2)\in J, |\theta\sigma^2 - \tilde{\theta}\tilde{\sigma}^2| \geq \varepsilon\}} \{l_{n,\text{tap}}(\tilde{\theta},\tilde{\sigma}^2) - l_{n,\text{tap}}(\theta,\sigma^2)\} \longrightarrow \infty\Big) = 1,$$

where $(\tilde{\theta},\tilde{\sigma}^2) \in J$ can be any fixed vector such that $\tilde{\theta}\tilde{\sigma}^2 = \theta_0\sigma_0^2$.

Ying (1991) has shown (A.37) for the log likelihood function $l_n(\theta,\sigma^2)$. More specifically, Ying (1991) showed that, uniformly in $(\theta,\sigma^2) \in J$ and $|\theta\sigma^2 - \tilde{\theta}\tilde{\sigma}^2| \geq \varepsilon$, with probability 1,

$$l_n(\tilde{\theta},\tilde{\sigma}^2) - l_n(\theta,\sigma^2) \geq \eta n + O_u(n^{1/2+\alpha}) \qquad \text{for any } \alpha > 0$$

[see the proof of Theorem 1 in Ying (1991), page 289]. Then, (A.37) follows because of (2.1) in Theorem 1.

Similarly, we can show (2.4) by using (2.2) and some asymptotic results in Ying (1991). We can write [see (3.10) and (3.11) in Ying (1991), page 291]

$$\frac{\partial}{\partial \theta}l_n(\theta,\sigma^2) = \frac{\sigma_0^2\theta_0}{\sigma^2\theta^2}\sum_{k=2}^{n} W_{k,n}^2 - \frac{n}{\theta} + O_u(1),$$

where

$$W_{k,n} = \frac{X(t_k) - e^{-\theta_0 \Delta_k}X(t_{k-1})}{\sigma_0\sqrt{1 - e^{-2\theta_0 \Delta_k}}}.$$

Note that $W_{k,n}$ depends only on the true parameters and are i.i.d. $N(0,1)$ for $k = 1,\ldots,n$.

Then, for any $(\theta,\sigma^2) \in J$, we have

$$\sigma^2\theta^2\frac{\partial}{\partial \theta}l_{n,\text{tap}}(\theta,\sigma^2) = \sigma_0^2\theta_0\sum_{k=2}^{n}(W_{k,n}^2 - 1) - n(\sigma^2\theta - \sigma_0^2\theta_0) + o_u(n^{1/2})$$

by Theorem 1. In particular, for $(\theta, \sigma^2) = (\hat{\theta}_{n,\text{tap}}, \hat{\sigma}^2_{n,\text{tap}})$, the left-hand side is zero. Therefore, we obtain

$$0 = \theta_0 \sigma_0^2 \sum_{k=2}^n (W_{k,n}^2 - 1) - n(\hat{\theta}_{n,\text{tap}} \hat{\sigma}^2_{n,\text{tap}} - \theta_0 \sigma_0^2) + o_u(n^{1/2}).$$

Since $W_{k,n}^2$, $k = 1, \ldots, n$, are i.i.d. $\chi_1^2$, we have

$$\sqrt{n}(\hat{\theta}_{n,\text{tap}} \hat{\sigma}^2_{n,\text{tap}} - \theta_0 \sigma_0^2) \;=\; \theta_0 \sigma_0^2 n^{-1/2} \sum_{k=2}^n (W^0{}_{k,n}^2 - 1) + o_u(1)$$

$$\xrightarrow{d} N(0, 2(\theta_0 \sigma_0^2)^2).$$

The proof is complete. $\square$

## APPENDIX B: PROOFS FOR SECTION 3

We will employ some known properties of equivalent Gaussian measures and will refer to Ibragimov and Rozanov (1978) frequently. Two measures $P_j, j = 0, 1$ are equivalent if they are absolutely continuous with respect to each other. Let $X(\mathbf{t})$, $\mathbf{t} \in D$ be Gaussian stationary under the two equivalent measures $P_j$, where $D$ is a bounded subset in $\mathbb{R}^d$ for some $d \geq 1$. Let $\mathbf{X}_n = (X(\mathbf{t}_1), \ldots, X(\mathbf{t}_n))'$ denote the observations in $D$, and let $p_j(x_1, \ldots, x_n)$ denote the density function of $\mathbf{X}_n$ under measure $P_j$ for $j = 0, 1$. Then, the Radon–Nikodym derivative $\rho_n = p_1(\mathbf{X}_n)/p_0(\mathbf{X}_n)$ has a limit $\rho$ with $P_0$-probability 1. In addition,

(B.1)
$$P_0(0 < \rho < \infty) = 1, \qquad \lim_{n \to \infty} E_0(\log \rho_n) = E_0(\log \rho) \quad \text{and}$$
$$-\infty < E_0(\log \rho) < \infty.$$

We refer the readers to Section III.2.1 of Ibragimov and Rozanov (1978) for these results. It follows that the log-likelihood ratio and its expectation are all bounded. With $P_0$-probability 1,

(B.2)
$$\log \rho_n = l_{n,1}(\mathbf{X}_n) - l_{n,0}(\mathbf{X}_n)$$
$$= -\frac{1}{2} \log \frac{\det \mathbf{V}_{1,n}}{\det \mathbf{V}_{0,n}} - \frac{1}{2} \mathbf{X}'_n (\mathbf{V}_{1,n}^{-1} - \mathbf{V}_{0,n}^{-1}) \mathbf{X}_n = O(1),$$

(B.3) $\quad E_0(\log \rho_n) = -\dfrac{1}{2} \log \dfrac{\det \mathbf{V}_{1,n}}{\det \mathbf{V}_{0,n}} - \dfrac{1}{2} E_0(\mathbf{X}'_n (\mathbf{V}_{1,n}^{-1} - \mathbf{V}_{0,n}^{-1}) \mathbf{X}_n) = O(1).$

The difference of these two equations yields

(B.4) $\quad \mathbf{X}'_n (\mathbf{V}_{1,n}^{-1} - \mathbf{V}_{0,n}^{-1}) \mathbf{X}_n - E_0(\mathbf{X}'_n (\mathbf{V}_{1,n}^{-1} - \mathbf{V}_{0,n}^{-1}) \mathbf{X}_n) = O(1) \qquad \text{a.s.}$

Before we proceed with the proof of the main results in Section 3, we will establish the following lemmas. For two functions $a(x), b(x)$, we write $a(x) \asymp b(x), x \to \infty$ if $-\infty < \liminf_{x \to \infty} a(x)/b(x) \leq \limsup_{x \to \infty} a(x)/b(x) < \infty$.

LEMMA 4. *Let $f_1(\lambda)$ be the spectral density corresponding to isotropic Matérn covariogram $K(h; \sigma_1^2, \theta_1)$ and $\tilde{f}_1(\lambda)$ be the spectral density corresponding to the tapered covariance function $\tilde{K}(h; \sigma_1^2, \theta_1) = K(h; \sigma_1^2, \theta_1) K_{\text{tap}}(h)$. Under condition (A3), there exists $r > 1$ such that*

$$\text{(B.5)} \qquad \frac{\tilde{f}_1(\lambda) - f_1(\lambda)}{f_1(\lambda)} = O(|\lambda|^{-r}) \qquad \text{as } |\lambda| \to \infty.$$

PROOF. Using the fact that Fourier transform of product of two functions is the convolution of their Fourier transforms, we have

$$\text{(B.6)} \qquad \tilde{f}_1(\lambda) = \int_{\mathbb{R}} f_1(x) f_{\text{tap}}(\lambda - x)\, dx,$$

where $f_{\text{tap}}$ is the spectral density corresponding to $K_{\text{tap}}$. It is seen that $\tilde{f}_1(\lambda)/f_1(\lambda)$ does not depend on $\sigma_1^2$ so that we can assume without loss of generality that $\sigma_1^2 = 1$. It suffices to consider the case that $\lambda > 0$, because $\tilde{f}_1(\lambda)$ is symmetric about $\lambda = 0$. Using $\int_{\mathbb{R}} f_{\text{tap}}(\lambda - x)\, dx = 1$ and breaking down these integrals over intervals $(-\infty, \lambda - \lambda^k] \cup [\lambda + \lambda^k, +\infty)$ and $(\lambda - \lambda^k, \lambda + \lambda^k)$ for any $k \in (0, 1)$, we have

$$\frac{\tilde{f}_1(\lambda)}{f_1(\lambda)} - 1 = \frac{\int_{|\lambda - x| \geq \lambda^k} f_1(x) f_{\text{tap}}(\lambda - x)\, dx}{f_1(\lambda)} - \int_{|\lambda - x| \geq \lambda^k} f_{\text{tap}}(\lambda - x)\, dx$$

$$+ \frac{\int_{|\lambda - x| < \lambda^k} (f_1(x) - f_1(\lambda)) f_{\text{tap}}(\lambda - x)\, dx}{f_1(\lambda)}$$

$$= T_1 + T_2 + T_3.$$

By condition (A3), we have

$$|T_1| \leq \frac{M}{(1 + \lambda^{2k})^{\nu + 1/2 + \varepsilon}} \frac{1}{f_1(\lambda)} \int f_1(x)\, dx.$$

The Matérn spectral density has a closed form

$$\text{(B.7)} \qquad f_1(\lambda) = \frac{c \theta_1^{2\nu}}{(\theta_1^2 + |\lambda|^2)^{\nu + 1/2}} \qquad \text{for } c = \frac{\Gamma(\nu + 1/2)}{\Gamma(\nu) \pi^{1/2}}.$$

In addition, $\int f_1(\lambda)\, d\lambda = \sigma_1^2 = 1$ is the variance. Then,

$$\text{(B.8)} \qquad |T_1| \leq \frac{M(\theta_1^2 + \lambda^2)^{\nu + 1/2}}{c \theta_1^{2\nu} (1 + \lambda^{2k})^{\nu + 1/2 + \varepsilon}}.$$

Similarly,

$$|T_2| \le \frac{M}{(\nu+\varepsilon)\lambda^{2k(v+\varepsilon)}}. \tag{B.9}$$

Since $\varepsilon > 1/2$ and $v + \varepsilon > 1$ by condition (A3), we can choose $k$ to be sufficiently close to 1, so that both $T_1$ and $T_2$ are $O(\lambda^{-r})$ for some $r > 1$.

To bound $T_3$, write, for some $\xi$ between $\lambda$ and $x$,

$$f_1(x) - f_1(\lambda) = f_1'(\lambda)(x-\lambda) + f_1''(\xi)\frac{(x-\lambda)^2}{2}.$$

Then,

$$T_3 = \frac{1}{f_1(\lambda)}\left(f_1'(\lambda)\int_{|x-\lambda|<\lambda^k}(x-\lambda)f_{\text{tap}}(\lambda-x)\,dx \right.$$
$$\left. + \int_{|\lambda-x|<\lambda^k} f_1''(\xi)\frac{(x-\lambda)^2}{2}f_{\text{tap}}(\lambda-x)\,dx\right).$$

The first term is 0 because the integrand is odd. For the second term, note

$$0 < f_1''(\xi) = \frac{c\theta_1^{2\nu}(2\nu+1)}{(\theta_1^2+\xi^2)^{\nu+3/2}}\left(\frac{(2\nu+3)\xi^2}{\theta_1^2+\xi^2}-1\right)$$
$$\le \frac{c\theta_1^{2\nu}(2\nu+3)^2}{(\theta_1^2+\xi^2)^{\nu+3/2}}.$$

Therefore, if $\lambda$ is sufficiently large, for $\xi$ lying between $x$ and $\lambda$, where $x$ is in the interval $|x-\lambda| < \lambda^k$,

$$0 < f_1''(\xi) \le \frac{c\theta_1^{2\nu}(2\nu+3)^2}{(\theta_1^2+(\lambda-\lambda^k)^2)^{\nu+3/2}} \le \frac{2c\theta_1^{2\nu}(2\nu+3)^2}{\lambda^{2\nu+3}}.$$

Then,

$$0 < T_3 < \frac{(2\nu+3)^2(\theta_1^2+\lambda^2)^{\nu+1/2}}{\lambda^{2\nu+3}}\int (x-\lambda)^2 f_{\text{tap}}(x-\lambda)\,dx.$$

Condition (A3) implies that $x^2 f_{\text{tap}}(x)$ is integrable. Then, $T_3 = O(\lambda^{-2})$. The proof is complete. □

LEMMA 5. *For any real number $r > 0$, there exists $\xi_r(\lambda)$ such that*

$$\xi_r(\lambda) = \int c_r(t)\exp(-i\lambda t)\,dt, \qquad 0 < |\xi_r(\lambda)|^2 \asymp |\lambda|^{-r}, \ \lambda \to \infty, \tag{B.10}$$

*where $c_r(t)$ is square integrable and has a compact support.*

PROOF. We only need to show the case $0 < r \leq 1$, because the product of any functions of the type given by (B.10) belongs to this type, due to the fact that the Fourier transform of convolution coincides with the product of Fourier transforms. Let

$$\xi_r(\lambda) = \int_{-1}^{1} e^{i\lambda t} |t|^{r/2-1} \, dt = 2 \int_0^1 \cos(\lambda t) t^{r/2-1} \, dt.$$

We will show that $\xi_r(\lambda)$ satisfies (B.10). We only need to prove it for $\lambda \geq 0$, because $\xi_r(\lambda)$ is symmetric about $\lambda = 0$ and $\xi_r(0) > 0$. Let $u = \lambda t$. We can write

$$\xi_r(\lambda) = 2\lambda^{-r/2} \int_0^{\lambda} \cos(u) u^{r/2-1} \, du.$$

Then, $\xi_r(\lambda)^2 \asymp |\lambda|^{-r}$ as $\lambda \to +\infty$, because $\cos(u)u^{r/2-1}$ is integrable for $0 < r \leq 1$.

Next, we will show $\xi_r(\lambda) > 0$ for any $\lambda > 0$. It suffices to show

(B.11)
$$y(\lambda) = \int_0^{\lambda} \cos(u) u^{-\delta} \, du > 0,$$

where $\delta = 1 - r/2 \in [1/2, 1)$. Note that $y'(\lambda) = \cos(\lambda)\lambda^{-\delta}$ and $y''(\lambda) = -\sin(\lambda) \times \lambda^{-\delta} - \delta \cos(\lambda)\lambda^{-\delta-1}$. Therefore, the minimum points are $\{2k\pi + 3\pi/2, k = 0, 1, \ldots\}$. So, we only need to show $y(2k\pi + 3\pi/2) > 0$, $k = 0, 1, \ldots$ by induction. First, using monotonicity of $\cos(u)$, we have

(B.12)
$$y\left(\frac{3\pi}{2}\right) = \int_0^{\pi/4} \cos(u)u^{-\delta} \, du + \int_{\pi/4}^{\pi/2} \cos(u)u^{-\delta} \, du$$
$$+ \int_{\pi/2}^{\pi} \cos(u)u^{-\delta} \, du + \int_{\pi}^{3\pi/2} \cos(u)u^{-\delta} \, du$$
$$\geq \cos(\pi/4) \frac{(\pi/4)^{1-\delta}}{1-\delta} + \frac{1}{(\pi/2)^{\delta}}\left(1 - \frac{\sqrt{2}}{2}\right) - \frac{1}{(\pi/2)^{\delta}} - \frac{1}{\pi^{\delta}}$$
$$\geq \frac{\sqrt{2}}{2}\sqrt{\pi} + \frac{2}{\pi}\left(1 - \frac{\sqrt{2}}{2}\right) - \sqrt{\frac{2}{\pi}} - \frac{1}{\sqrt{\pi}} > 0.$$

Next, suppose $y(2(k-1)\pi + 3\pi/2) > 0$, for $k \geq 1$, then

$$y(2k\pi + 3\pi/2) = y(2(k-1)\pi + 3\pi/2)$$
$$+ \int_{2k\pi - \pi/2}^{2k\pi + \pi/2} \cos(u)u^{-\delta} \, du + \int_{2k\pi + \pi/2}^{2k\pi + 3\pi/2} \cos(u)u^{-\delta} \, du$$
$$= y(2(k-1)\pi + 3\pi/2) + \int_{2k\pi - \pi/2}^{2k\pi + \pi/2} \cos(u)u^{-\delta} \, du$$

$$-\int_{2k\pi-\pi/2}^{2k\pi+\pi/2} \cos(u)(u+\pi)^{-\delta}\,du$$

$$= y(2(k-1)\pi + 3\pi/2)$$

$$+ \int_{2k\pi-\pi/2}^{2k\pi+\pi/2} \cos(u)(u^{-\delta} - (u+\pi)^{-\delta})\,du.$$

The integral is positive because the integrand is positive. This completes the proof of Lemma 5. □

PROOF OF THEOREM 3. Write the Cholesky decomposition of $\mathbf{V}_{0,n} = \mathbf{L}\mathbf{L}'$ for some lower triangular matrix $\mathbf{L}$. Let $\mathbf{Q}$ be an orthogonal matrix such that

$$\mathbf{Q}\mathbf{L}^{-1}\mathbf{V}_{1,n}\mathbf{L}'^{-1}\mathbf{Q}' = \operatorname{diag}\{\sigma_{1,n}^2,\ldots,\sigma_{n,n}^2\}.$$

Then,

$$\mathbf{Q}'\mathbf{L}'\mathbf{V}_{1,n}^{-1}\mathbf{L}\mathbf{Q} = \operatorname{diag}\{1/\sigma_{1,n}^2,\ldots,1/\sigma_{n,n}^2\}.$$

Taking the trace of both sides, we have

$$\operatorname{trace}(\mathbf{V}_{0,n}\mathbf{V}_{1,n}^{-1}) = \sum_{i=1}^{n} 1/\sigma_{k,n}^2.$$

Hence,

$$E_0(\mathbf{X}_n'(\mathbf{V}_{1,n}^{-1} - \mathbf{V}_{0,n}^{-1})\mathbf{X}_n) = \operatorname{trace}(\mathbf{V}_{0,n}\mathbf{V}_{1,n}^{-1}) - n = \sum_{k=1}^{n}\left(\frac{1}{\sigma_{k,n}^2} - 1\right).$$

Let $\mathbf{e}_n = \mathbf{Q}\mathbf{L}^{-1}\mathbf{X}_n$. Obviously,

(B.13) $\qquad E_0\mathbf{e}_n\mathbf{e}_n' = \mathbf{I}_n, \qquad E_1\mathbf{e}_n\mathbf{e}_n' = \operatorname{diag}\{\sigma_{1,n}^2,\ldots,\sigma_{n,n}^2\}.$

Equation (3.7) follows if, for any orthogonal sequence $\{\eta_k, k=1,2,\ldots\}$ in the Hilbert space $L_D^2(dP_0)$ spanned by $X(t), t \in D$ under the covariance inner product corresponding to $P_0$, there exists a constant $M > 0$ independent of $\eta_k, k = 1, 2, \ldots$, such that

(B.14) $$\sum_{k=1}^{\infty}\left|\frac{1}{E_1\eta_k^2} - 1\right| < M.$$

One important technique to prove (B.14) is to write, for any $s, t \in D$,

(B.15) $\quad E_1 X(t)X(s) - E_0 X(t)X(s) = \int_{\mathbb{R}}\int_{\mathbb{R}} e^{i(\lambda s - \mu t)}\Phi(\lambda,\mu)\,d\lambda\,d\mu,$

where $\Phi(\lambda,\mu)$ is square integrable with respect to Lebesgue measure on $\mathbb{R}^2$. For any bounded region $D$, the existence of such a function $\Phi$ and, therefore, the equivalence of $P_0$ and $P_1$, are shown in Ibragimov and Rozanov (1978), page 104, Theorem 17, under the assumption that the function $h(\lambda)$ in (3.6) is square integrable. However, we will show, under the assumption of this lemma that $h(\lambda)$ is integrable, $\Phi$ takes a particular form

$$\Phi(\lambda,\mu) = \overline{\Phi_1(\lambda)}\Phi_2(\mu) \int_T e^{-i(\lambda-\mu)\omega}\,d\omega \tag{B.16}$$

for some functions $\Phi_j(\lambda), \lambda \in \mathbb{R}$ such that $\int |\Phi_j(\lambda)|^2/f_0(\lambda)\,d\lambda < \infty$, $j=1,2$, and a compact interval $T$ that is solely determined by $r_1$ and $r_2$. This particular form is central to the proof, and we will establish it at the end of this proof. We now proceed by assuming it is true.

Let $dZ_0(\lambda)$ denote the stochastic orthogonal measure so that $X(t)$ has the spectral representation under measure $P_0$; that is, $X(t) = \int \exp(-i\lambda t)\,dZ_0(\lambda)$. Then, for any $\eta \in L^2_D(dP_0)$, there is a function $\phi(\lambda)$ such that $\eta = \int \phi(\lambda)\,dZ_0(\lambda)$ and $E_0\eta^2 = \int |\phi(\lambda)|^2 f_0(\lambda)\,d\lambda$. We first show

$$E_1\eta^2 = \int |\phi(\lambda)|^2 f_1(\lambda)\,d\lambda, \tag{B.17}$$

$$E_1\eta^2 - E_0\eta^2 = \int\int \overline{\phi(\lambda)}\phi(\mu)\Phi(\lambda,\mu)\,d\lambda\,d\mu. \tag{B.18}$$

Indeed, the two equations hold for $\eta = X(t) = \int \exp(-i\lambda t)\,dZ_0(\lambda)$ for any $t \in D$ [assuming (B.15) is true]. Consequently, they hold for any linear combination

$$\eta = \sum_{j=1}^{J} c_j X(t_j) = \int \phi(\lambda)\,dZ_0(\lambda)$$

for any $J$ and $t_1,\ldots,t_J \in D$, where $\phi(\lambda) = \sum_{j=1}^{J} c_j e^{-i\lambda t_j}$.

For any $\eta \in L^2_D(dP_0)$, we can find a sequence of finite linear combinations of $X(t), t \in D$, say, $\eta_m, m=1,2,\ldots$, such that $\lim_{m\to\infty} E_0(\eta-\eta_m)^2 = 0$. If $\eta_m = \int \phi_m(\lambda)\,dZ_0(\lambda)$, we have

$$E_0(\eta-\eta_m)^2 = \int |\phi(\lambda) - \phi_m(\lambda)|^2 f_0(\lambda)\,d\lambda \to 0. \tag{B.19}$$

Then,

$$\int |\phi(\lambda) - \phi_m(\lambda)|^2 f_1(\lambda)\,d\lambda = \int |\phi(\lambda) - \phi_m(\lambda)|^2 f_0(\lambda)(1+h(\lambda))\,d\lambda \to 0,$$

because $h = (f_1 - f_0)/f_0$ is bounded. It follows that $\eta_m$ converges in $L^2(dP_1)$ norm to some variable $\tilde{\eta}$ because $E_1(\eta_l - \eta_m)^2 = \int |\phi_l(\lambda) - \phi_m(\lambda)|^2 f_1(\lambda)\,d\lambda \to 0$ as $l,m \to \infty$. Then,

$$E_1\tilde{\eta}^2 = \lim_{m\to\infty} E_1\eta_m^2 = \lim_{m\to\infty} \int |\phi_m(\lambda)|^2 f_1(\lambda)\,d\lambda = \int |\phi(\lambda)|^2 f_1(\lambda)\,d\lambda.$$

Since $L_2$ convergence implies convergence in probability, we have $\eta_m \to \tilde{\eta}$ in probability $P_1$. On the other hand, $\eta_m \to \eta$ in probability $P_0$ and, consequently, in probability $P_1$, due to the equivalence of the two probabilities. Then, we must have $P_1(\eta = \tilde{\eta}) = 1$ and $E_1 \eta^2 = E_1 \tilde{\eta}^2$. We have proved (B.17). To show (B.18), note that

$$\left| \iint \overline{\phi_m(\lambda)} \phi_m(\mu) \Phi(\lambda, \mu) \, d\lambda \, d\mu - \iint \overline{\phi(\lambda)} \phi(\mu) \Phi(\lambda, \mu) \, d\lambda \, d\mu \right|$$

(B.20)
$$\leq \iint |(\overline{\phi_m(\lambda)} - \overline{\phi(\lambda)}) \phi_m(\mu)| |\Phi(\lambda, \mu)| \, d\lambda \, d\mu$$
$$+ \iint |(\phi_m(\mu) - \phi(\mu)) \overline{\phi(\lambda)}| |\Phi(\lambda, \mu)| \, d\lambda \, d\mu,$$

where the first term tends to zero, because Cauchy–Schwarz inequality implies its square is bounded by

$$|T|^2 \int |\overline{\phi_m(\lambda)} - \overline{\phi(\lambda)}|^2 f_0(\lambda) \, d\lambda \int |\phi_m(\mu)|^2 f_0(\mu) \, d\mu$$
$$\times \int \frac{|\Phi_1(\lambda)|^2}{f_0(\lambda)} \, d\lambda \int \frac{|\Phi_2(\mu)|^2}{f_0(\mu)} \, d\mu \to 0$$

by (B.19) and square integrability of $\Phi_1(\lambda)/\sqrt{f_0(\lambda)}$, where and hereafter $|T|$ stands for the length of finite interval $T$. Similarly, we can show the second term in (B.20) also tends to zero. Therefore, (B.18) is now proved by taking the limit of $E_1 \eta_m^2 - E_0 \eta_m^2$ and $\iint \overline{\phi_m(\lambda)} \phi_m(\mu) \Phi(\lambda, \mu) \, d\lambda \, d\mu$.

Applying (B.18) to the orthonormal sequence $\eta_k = \int \phi_k(\lambda) \, dZ_0(\lambda)$, $k = 1, 2, \ldots$, we have

$$E_1 \eta_k^2 - 1 = \iint \overline{\phi_k(\lambda)} \phi_k(\mu) \int_T e^{-i(\lambda - \mu)\omega} \, d\omega \, \overline{\Phi_1(\lambda)} \Phi_2(\mu) \, d\lambda \, d\mu$$
$$= \int_T \overline{A_{1,k}(\omega)} A_{2,k}(\omega) \, d\omega,$$

where $A_{j,k}(\omega) = \int \phi_k(\lambda) \exp(i\lambda\omega) \Phi_j(\lambda) \, d\lambda$. Since (3.6) and continuity of $f_j$ imply that $f_1(\lambda) > C f_0(\lambda)$ for some constant $C > 0$, we have $E_1 \eta_k^2 > C E_0 \eta_k^2 = C$ and, therefore,

$$|1/E_1 \eta_k^2 - 1| \leq |E_1 \eta_k^2 - 1|/C \leq (1/2C) \sum_{j=1}^{2} \int_T |A_{j,k}(\omega)|^2 \, d\omega.$$

In view that $A_{j,k}(\omega)$ is the inner product of the two integrable functions $\phi_k(\lambda) f_0(\lambda)^{1/2}$ and $\exp(i\lambda\omega) \Phi_j(\lambda)/f_0(\lambda)^{1/2}$ in $L^2(d\lambda)$, and that $\phi_k(\lambda) f_0(\lambda)^{1/2}$, $k = 1, 2, \ldots$, is an orthonormal sequence in $L^2(d\lambda)$ [because $E_0 \eta_l \eta_k = \int \phi_l(\lambda) \times \phi_k(\lambda) f_0(\lambda) \, d\lambda$], we have, by Bessel's inequality,

$$\sum_{k=1}^{\infty} |A_{j,k}(\omega)|^2 \leq \int |\Phi_j(\lambda)|^2 / f_0(\lambda) \, d\lambda < \infty.$$

It follows that
$$\sum_{k=1}^{\infty} |E_1\eta_k^2 - 1| \le (1/2C) \sum_{j=1}^{2} \int_T \sum_{k=1}^{\infty} |A_{j,k}(\omega)|^2 \, d\omega$$
$$\le (|T|/2C) \sum_{j=1}^{2} \int |\Phi_j(\lambda)|^2/f_0(\lambda) \, d\lambda < \infty.$$

We just need to show (B.15) and (B.16) to complete the proof. We will employ the following well-known properties of Fourier transform. For any square integrable functions (with respect to Lebesgue measure) $\varphi_j(\boldsymbol{\lambda})$, $\boldsymbol{\lambda} \in \mathbb{R}^d$, there are square integrable functions $a_j(\mathbf{t}), \mathbf{t} \in \mathbb{R}^d$ such that
$$\varphi_j(\boldsymbol{\lambda}) = \int_{\mathbb{R}^d} \exp(-i\boldsymbol{\lambda}'\mathbf{t})a_j(\mathbf{t}) \, d\mathbf{t}, \qquad j=1,2.$$
Furthermore,

(B.21) $$\varphi_1(\boldsymbol{\lambda})\varphi_2(\boldsymbol{\lambda}) = \int_{\mathbb{R}^d} \exp(-i\boldsymbol{\lambda}'\mathbf{t})(a_1 * a_2)(\mathbf{t}) \, d\mathbf{t},$$

(B.22) $$\int_{\mathbb{R}^d} \exp(i\boldsymbol{\lambda}'\mathbf{t})\varphi_1(\boldsymbol{\lambda})\varphi_2(\boldsymbol{\lambda}) \, d\boldsymbol{\lambda} = (2\pi)^d (a_1 * a_2)(\mathbf{t}),$$

where all the equalities are in the $L^2(d\boldsymbol{\lambda})$ sense, and $a_1 * a_2$ is the convolution; that is,
$$a_1 * a_2(\mathbf{t}) = \int_{\mathbb{R}^d} a_1(\mathbf{s})a_2(\mathbf{t}-\mathbf{s}) \, d\mathbf{s}.$$

By Lemma 5, there exists a continuous and square integrable function $\xi_j(\lambda)$ ($j = 1,2$) such that

(B.23) $$\xi_j(\lambda) = \int c_j(t) \exp(-i\lambda t) \, dt, \qquad 0 < |\xi_j(\lambda)|^2 \asymp |\lambda|^{-r_j}, \; \lambda \to \infty,$$

for some square integrable function $c_j(t)$ that has a compact support [i.e., $c_j(t)$ is 0 outside a compact set].

Let $\xi(\lambda) = (f_0(\lambda) - f_1(\lambda))/|\xi_1(\lambda)|^2$. Then, $\xi(\lambda)$ is square integrable by the assumption of the theorem and the properties of $\xi_1(\lambda)$. Therefore, we can write, for some square integrable function $c(t)$, $\xi(\lambda) = \int \exp(-i\lambda t)c(t) \, dt$. Furthermore, for all $s, t$,

(B.24)
$$E_0 X(s)X(t) - E_1 X(s)X(t) = \int e^{i\lambda(s-t)}(f_0(\lambda) - f_1(\lambda)) \, d\lambda$$
$$= \int e^{i\lambda(s-t)}\xi(\lambda)|\xi_1(\lambda)|^2 \, d\lambda,$$

which we will denote by $b(s,t)$. By (B.21),
$$|\xi_1(\lambda)|^2 = \int \exp(-i\lambda t)\left(\int c_1(z)\overline{c_1(z-t)} \, dz\right) dt.$$

Applying (B.22) to $\xi(\lambda)$ and $|\xi_1(\lambda)|^2$, we get

$$b(s,t) = 2\pi \int_{\mathbb{R}} c(w) \int_{\mathbb{R}} c_1(z)\overline{c_1(-(s-t-w-z))} \, dz \, dw$$
$$(B.25)$$
$$= 2\pi \int_{\mathbb{R}^2} c(u-v) c_1(s-u) \overline{c_1(t-v)} \, du \, dv,$$

which holds for all $s, t \in \mathbb{R}$. If we restrict $s, t$ to the compact set $D$, the integral (B.25) is an integral over a compact set, say, $\Delta \times \Delta$. This is because $c_1$ is 0 outside a compact interval.

Next, we write $c(t)$ as a convolution of two functions. For this purpose, we write $\xi(\lambda) = \xi_2(\lambda)\xi_3(\lambda)$. Then, $\xi_3(\lambda)$ so defined is square integrable from assumptions and (B.23) and, therefore, can be written as

$$\xi_3(\lambda) = \int \exp(-i\lambda t) c_3(t) \, dt.$$

Then, $c = c_2 * c_3$ and, consequently,

$$(B.26) \quad c(u-v) = \int c_2(x) c_3(u-v-x) \, dx = \int c_2(u-\omega) c_3(\omega-v) \, d\omega.$$

Since we are only interested in $b(s,t)$ for $s, t \in D$ and, consequently, only interested in $c(u-v)$ for $u, v \in \Delta$, we will restrict both $u, v$ to the interval $\Delta$, so that the second interval in (B.26) is an integral on a finite interval, say, $T$, because $c_2$ has a compact support. Define the bivariate function

$$a(u,v) = \int_T c_2(u-\omega) c_3(\omega-v) \, d\omega, \qquad u, v \in \mathbb{R},$$

which is square integrable because

$$|a(u,v)|^2 \leq |T| \int_T |c_2(u-\omega)|^2 |c_3(\omega-v)|^2 \, d\omega$$

and both $c_2$ and $c_3$ are square integrable. In addition, for $u, v \in \Delta$, we have, from (B.26),

$$a(u,v) = c(u-v).$$

We therefore have shown that, for $s, t \in D$,

$$b(s,t) = 2\pi \int_{\mathbb{R}^2} a(u,v) c_1(s-u) \overline{c_1(t-v)} \, du \, dv.$$

Note that the integral is a convolution of functions of $(u,v)$. Applying (B.22), we get

$$2\pi b(s,t) = \int \exp(i(\lambda s + \mu t)) \varphi_1(\lambda, \mu) \varphi_2(\lambda, \mu) \, d\lambda \, d\mu,$$

where
$$\varphi_1(\lambda,\mu) = \int_{\mathbb{R}^2} a(u,v) e^{-i(u\lambda+v\mu)}\,du\,dv,$$
$$\varphi_2(\lambda,\mu) = \int_{\mathbb{R}^2} c_1(u)\overline{c_1(v)} e^{-i(u\lambda+v\mu)}\,du\,dv.$$

Clearly,
$$\varphi_2(\lambda,\mu) = \xi_1(\lambda)\overline{\xi_1(-\mu)}.$$

Now,
$$\begin{aligned}
\varphi_1(\lambda,\mu) &= \int_{\mathbb{R}^2} a(u,v) e^{-i(u\lambda+v\mu)}\,du\,dv \\
&= \int_T \int_{\mathbb{R}^2} c_2(u-\omega) c_3(\omega-v) e^{-i(u\lambda+v\mu)}\,du\,dv\,d\omega \\
&= \int_T \int_{\mathbb{R}^2} c_2(x) c_3(-y) e^{-i((x+\omega)\lambda+(y+\omega)\mu)}\,dx\,dy\,d\omega \\
&= \int_T \left( \int_{\mathbb{R}^2} c_2(x) e^{-ix\lambda} c_3(-y) e^{-iy\mu}\,dx\,dy \right) e^{-i(\lambda+\mu)\omega}\,d\omega \\
&= \xi_2(\lambda)\xi_3(-\mu) \int_T e^{-i(\lambda+\mu)\omega}\,d\omega.
\end{aligned}$$

Hence,
$$\begin{aligned}
b(s,t) &= \frac{1}{2\pi} \int_{\mathbb{R}^2} e^{i(\lambda s+\mu t)} \xi_1(\lambda)\xi_2(\lambda)\overline{\xi_1(-\mu)}\xi_3(-\mu) \int_T e^{-i(\lambda+\mu)\omega}\,d\omega\,d\lambda\,d\mu \\
&= \frac{1}{2\pi} \int_{\mathbb{R}^2} \exp(i(\lambda s - \mu t))\overline{\Phi_1(\lambda)}\Phi_2(\mu) \int_T \exp(-i(\lambda-\mu)\omega)\,d\omega\,d\lambda\,d\mu
\end{aligned}$$

for $\Phi_1(\lambda) = \overline{\xi_1(\lambda)\xi_2(\lambda)}$ and $\Phi_2(\mu) = \overline{\xi_1(\mu)}\xi_3(\mu)$. Clearly, $\int |\Phi_j(\lambda)|^2/f_0(\lambda)\,d\lambda < \infty$ by the assumption of the theorem, (B.23) and the square-integrability of $\xi_2$ and $\xi_3$. The proof is complete. $\square$

PROOF OF THEOREM 4. Let $\sigma_1^2$ be such that $\sigma_1^2\theta_1^{2\nu} = \sigma_0^2\theta_0^{2\nu}$, and let $P_j$ be the probability measure under which the process has a Matérn covariogram with parameters $(\theta_j, \sigma_j^2)$ for $j = 0, 1$. Then, $P_0 \equiv P_1$ by Theorem 2 in Zhang (2004). Consequently, we only need to show that (3.8) and (3.9) hold, almost surely, with respect to $P_1$.

Let $f_1(\lambda) = f_1(\lambda; \theta_1, \sigma_1^2)$ be the spectral density under measure $P_1$ and $f_2(\lambda)$ the corresponding tapered spectral density as defined in Lemma 4, from which we see that, for some constant $c > 0$,
$$\int_{|\lambda|>c} \left| \frac{f_2(\lambda) - f_1(\lambda)}{f_1(\lambda)} \right|^2 d\lambda < \infty,$$

which is a sufficient condition for the equivalence of the measures $P_1$ and $P_2$ where $P_2$ is the measure corresponding the tapered spectral density $f_2$.

Let $\mathbf{V}_{j,n}$, $j=1,2$, be the covariance matrix corresponding to the spectral densities $f_j$ that depend on $\sigma_1^2$ and $\theta_1$ and do not depend on $\sigma^2$. For any $\sigma^2$, we have

$$
\text{(B.27)} \quad \begin{aligned} l_{n,\text{tap}}(\theta_1,\sigma^2) &- l_n(\theta_1,\sigma^2) \\ &= -\log(\det \mathbf{V}_{2,n}/\det \mathbf{V}_{1,n}) - \frac{\sigma_1^2}{\sigma^2}\mathbf{X}_n'(\mathbf{V}_{2,n}^{-1} - \mathbf{V}_{1,n}^{-1})\mathbf{X}_n. \end{aligned}
$$

Split it into three additive terms as follows:

$$
\text{(B.28)} \quad \begin{aligned} &\left[-\log\frac{\det \mathbf{V}_{2,n}}{\det \mathbf{V}_{1,n}} - E_1(\mathbf{X}_n'(\mathbf{V}_{2,n}^{-1} - \mathbf{V}_{1,n}^{-1})\mathbf{X}_n)\right] \\ &+ \left(1 - \frac{\sigma_1^2}{\sigma^2}\right)E_1(\mathbf{X}_n'(\mathbf{V}_{2,n}^{-1} - \mathbf{V}_{1,n}^{-1})\mathbf{X}_n) \end{aligned}
$$

$$
\text{(B.29)} \quad \begin{aligned} &- \frac{\sigma_1^2}{\sigma^2}[\mathbf{X}_n'(\mathbf{V}_{2,n}^{-1} - \mathbf{V}_{1,n}^{-1})\mathbf{X}_n - E_1(\mathbf{X}_n'(\mathbf{V}_{2,n}^{-1} - \mathbf{V}_{1,n}^{-1})\mathbf{X}_n)] \\ &= I_1 + I_2 - I_3. \end{aligned}
$$

Because $P_1 \equiv P_2$, the first term is bounded as we discussed previously in (B.3). Similarly, by (B.2), the third term $I_3$ is bounded uniformly in $\sigma^2 \in [w,v]$, almost surely. The second term $I_2$ is also bounded uniformly in $\sigma^2 \in [w,v]$ because, by Theorem 3,

$$
\text{(B.30)} \quad E_1(\mathbf{X}_n'(\mathbf{V}_{2,n}^{-1} - \mathbf{V}_{1,n}^{-1})\mathbf{X}_n) = O(1).
$$

Therefore (3.8) is proved. To show (3.9), first observe

$$
\text{(B.31)} \quad \frac{\partial}{\partial \sigma^2}l_{n,\text{tap}}(\theta_1,\sigma^2) - \frac{\partial}{\partial \sigma^2}l_n(\theta_1,\sigma^2) = -\frac{\sigma_1^2}{\sigma^4}(\mathbf{X}_n'\mathbf{V}_{2,n}^{-1}\mathbf{X}_n - \mathbf{X}_n'\mathbf{V}_{1,n}^{-1}\mathbf{X}_n),
$$

which can be rewritten as

$$
\text{(B.32)} \quad \begin{aligned} &\frac{\sigma_1^2}{\sigma^4}[\mathbf{X}_n'(\mathbf{V}_{2,n}^{-1} - \mathbf{V}_{1,n}^{-1})\mathbf{X}_n - E_1(\mathbf{X}_n'(\mathbf{V}_{2,n}^{-1} - \mathbf{V}_{1,n}^{-1})\mathbf{X}_n)] \\ &- \frac{\sigma_1^2}{\sigma^4}E_1(\mathbf{X}_n'(\mathbf{V}_{2,n}^{-1} - \mathbf{V}_{1,n}^{-1})\mathbf{X}_n). \end{aligned}
$$

Then, (3.9) immediately follows (B.4) and Theorem 3. □

PROOF OF THEOREM 5. As $\sigma_1^2 = \sigma_0^2(\theta_0/\theta_1)^{2\nu}$, we only need to show

$$
\text{(B.33)} \quad \sqrt{n}\left(\frac{\hat{\sigma}_n^2}{\sigma_1^2} - 1\right) \xrightarrow{d} N(0,2),
$$

$$
\text{(B.34)} \quad \sqrt{n}\left(\frac{\hat{\sigma}_{n,\text{tap}}^2}{\sigma_1^2} - 1\right) \xrightarrow{d} N(0,2).
$$

Let $\mathbf{V}_{j,n}$ be the covariance matrix of $\mathbf{X}_n$ corresponding to parameter values $(\theta_i, \sigma_i^2)$, $j = 0, 1$. Write $\mathbf{V}_{1,n} = \sigma_1^2 \mathbf{R}_{1,n}$, where $\mathbf{R}_{1,n}$ is the correlation matrix. First, we note that $\hat{\sigma}_n^2$ has a closed form express

$$\hat{\sigma}_n^2 = \frac{1}{n}\mathbf{X}_n' \mathbf{R}_{1,n}^{-1} \mathbf{X}_n \tag{B.35}$$

that can be derived straightforwardly from the maximization. Then,

$$\sqrt{n}\left(\frac{\hat{\sigma}_n^2}{\sigma_1^2} - 1\right) = \sqrt{n}\left(\frac{\mathbf{X}_n' \mathbf{R}_{1,n}^{-1} \mathbf{X}_n}{\sigma_1^2 n} - 1\right)$$
$$= (1/\sqrt{n})(\mathbf{X}_n'(\mathbf{V}_{1,n}^{-1} - \mathbf{V}_{0,n}^{-1})\mathbf{X}_n) + \sqrt{n}\left(\frac{\mathbf{X}_n' \mathbf{V}_{0,n}^{-1} \mathbf{X}_n}{n} - 1\right). \tag{B.36}$$

Since $\mathbf{V}_{0,n}^{-1/2} \mathbf{X}_n$ consists of i.i.d. $N(0,1)$ variables, $\mathbf{X}_n' \mathbf{V}_{0,n}^{-1} \mathbf{X}_n$ is the sum of i.i.d. variables having a $\chi_1^2$ distribution. The central limit theorem implies that the second term in (B.36) converges in distribution to $N(0, 2)$.

Equation (3.10) in Theorem 5 follows if the first term is shown to be bounded almost surely with respect to $P_0$. In view of (B.4), if suffices to show that

$$E_0(\mathbf{X}_n'(\mathbf{V}_{1,n}^{-1} - \mathbf{V}_{0,n}^{-1})\mathbf{X}_n) = O(1).$$

To this end, we only need to verify that conditions of Theorem 3 are satisfied. The Matérn spectral density (B.7) satisfies, as $\lambda \to \infty$,

$$0 < f(\lambda; \theta_i, \sigma_i^2) \sim |\lambda|^{-(2\nu+1)}. \tag{B.37}$$

Moreover, in view of $\sigma_0^2 \theta_0^{2\nu} = \sigma_1^2 \theta_1^{2\nu}$,

$$h(\lambda) = \frac{f_1(\lambda)}{f_0(\lambda)} - 1 = \left(\frac{\theta_0^2 + \lambda^2}{\theta_1^2 + \lambda^2}\right)^{\nu+d/2} - 1 = \left(1 + \frac{\theta_0^2 - \theta_1^2}{\theta_1^2 + \lambda^2}\right)^{\nu+d/2} - 1,$$
(B.38)
where $f_i(\lambda)$ stands for $f(\lambda; \theta_i, \sigma_i^2)$, $i = 0, 1$. Using the Taylor expansion, we can get

$$h(\lambda) \sim |\lambda|^{-2}.$$

Hence, (3.10) is proved.

Next, we derive the asymptotic distribution of the tapered MLE $\hat{\sigma}_{n,\text{tap}}^2$. Similar to $\hat{\sigma}_n^2$, the tapered MLE $\hat{\sigma}_{n,\text{tap}}^2$ takes the closed form

$$\hat{\sigma}_{n,\text{tap}}^2 = \frac{1}{n}\mathbf{X}_n' \tilde{\mathbf{R}}_{1,n}^{-1} \mathbf{X}_n,$$

where $\tilde{\mathbf{R}}_{1,n}$ is the tapered correlation matrix corresponding to $\mathbf{R}_{1,n}$. It follows (3.10) in Theorem 4 that

$$\mathbf{X}_n' \tilde{\mathbf{R}}_{1,n}^{-1} \mathbf{X}_n - \mathbf{X}_n' \mathbf{R}_{1,n}^{-1} \mathbf{X}_n = O(1).$$

Then,

$$\hat{\sigma}^2_{n,\text{tap}} = \frac{1}{n}\mathbf{X}'_n \mathbf{R}^{-1}_{1,n} \mathbf{X}_n + O(1/n) = \hat{\sigma}^2_n + O(1/n).$$

It follows immediately that $\hat{\sigma}^2_{n,\text{tap}}$ and $\hat{\sigma}^2_n$ have the same asymptotic distribution. The proof is complete. $\square$

**Acknowledgments.** The authors thank the editors and two referees for their constructive comments.

## REFERENCES

ABRAMOWITZ, M. and STEGUN, I. (1967). *Handbook of Mathematical Functions*. Government Printing Office, Washington, DC.

BELLMAN, R. E. (1970). *Introduction to Matrix Analysis*. McGraw-Hill, New York. MR0258847

CHEN, H.-S., SIMPSON, D. G. and YING, Z. (2000). Fixed-domain asymptotics for a stochastic process model with measurement error. *Statist. Sinica* **10** 141–156. MR1742105

DAVIS, T. A. (2006). *Direct Methods for Sparse Linear Systems*. SIAM, Philadelphia. MR2270673

FURRER, R., GENTON, M. G. and NYCHKA, D. (2006). Covariance tapering for interpolation of large spatial datasets. *J. Comput. Graph. Statist.* **15** 502–523. MR2291261

GILBERT, J. R., MOLER, C. and SCHREIBER, R. (1992). Sparse matrices in MATLAB: Design and implementation. *SIAM J. Matrix Anal. Appl.* **13** 333–356. MR1146669

GNEITING, T. (1999). Radial positive definite functions generated by Euclids hat. *J. Multivariate Anal.* **69** 88–119. MR1701408

GNEITING, T. (2002). Compactly supported correlation functions. *J. Multivariate Anal.* **83** 493–508. MR1945966

GRAYBILL, F. A. (1983). *Matrices with Applications in Statistics*, 2nd ed. Wadsworth, Belmont, CA. MR0682581

IBRAGIMOV, I. A. and ROZANOV, Y. A. (1978). *Gaussian Random Processes*. Springer, New York. MR0543837

KAUFMAN, C., SCHERVISH, M. and NYCHKA, D. (2008). Covariance tapering for likelihood-based estimation in large spatial datasets. *J. Amer. Statist. Assoc.* **103** 1545–1555.

LEMONE, M. A., CHEN, F., ALFIERI, J. G., CUENCA, R., HAGIMOTO, Y., BLANKEN, P. D., NIYOGI, D., KANG, S., DAVIS, K. and GROSSMAN, R. (2007). NCAR/CU surface vegetation observation network during the international H2O project 2002 field campaign. *Bull. Amer. Meteor. Soc.* **88** 65–81.

LOH, W.-L. (2005). Fixed-domain asymptotics for a subclass of Matérn-type Gaussian random fields. *Ann. Statist.* **33** 2344–2394. MR2211089

MARDIA, K. V. and MARSHALL, R. J. (1984). Maximum likelihood estimation of models for residual covariance in spatial regression. *Biometrika* **71** 135–146. MR0738334

MIRSKY, L. (1955). *An Introduction to Linear Algebra*. Oxford Univ. Press, Oxford. MR0074364

PISSANETZKY, S. (1984). *Sparse Matrix Technology*. Academic Press, London. MR0751237

RIPLEY, B. D. (1981). *Spatial Statistics*. Wiley, New York. MR0624436

SEARLE, S. R. (1971). *Linear Models*. Wiley, New York. MR0293792

Stassberg, D., LeMone, M. A., Warner, T. and Alfieri, J. G. (2008). Comparison of observed 10 m wind speeds to those based on Monin–Obukhov similarity theory using aircraft and surface data from the International H2O project. *Mon. Wea. Rev.* **136** 964–972.

Stein, M. L. (1988). Asymptotically efficient prediction of a random field with a misspecified covariance function. *Ann. Statist.* **16** 55–63. MR0924856

Stein, M. L. (1990a). Uniform asymptotic optimality of linear predictions of a random field using an incorrect second-order structure *Ann. Statist.* **18** 850–872. MR1056340

Stein, M. L. (1990b). Bounds on the efficiency of linear predictions using an incorrect covariance function. *Ann. Statist.* **18** 1116–1138. MR1062701

Stein, M. L. (1990c). A comparison of generalized cross validation and modified maximum likelihood for estimating the parameters of a stochastic process. *Ann. Statist.* **18** 1139–1157. MR1062702

Stein, M. L. (1999a). Predicting random fields with increasing dense observations. *Ann. Statist.* **9** 242–273. MR1682572

Stein, M. L. (1999b). *Interpolation of Spatial Data: Some Theory for Kriging*. Springer, New York. MR1697409

Weckworth, T. M., Parsons, D. B., Koch, S. E., Moore, J. A., LeMone, M. A., Demoz, B. B., Flamant, C., Geerts, B., Wang, J. and Feltz, W. F. (2004). An overview of the international H2O project (IHOP_2002) and some preliminary highlights. *Bull. Amer. Meteor. Soc.* **85** 253–277.

Wendland, H. (1995). Piecewise polynomial, positive definite and compactly supported radial functions of minimal degree. *Adv. Comput. Math.* **4** 389–396. MR1366510

Wendland, H. (1998). Error estimates for interpolation by compactly supported radial basis functions of minimal degree. *J. Approx. Theory* **93** 258–272. MR1616781

Wu, Z. M. (1995). Compactly supported positive definite radial functions. *Adv. Comput. Math.* **4** 283–292. MR1357720

Ying, Z. (1991). Asymptotic properties of a maximum likelihood estimator with data from a Gaussian process. *J. Multivariate Anal.* **36** 280–296. MR1096671

Ying, Z. (1993). Maximum likelihood estimation of parameters under a spatial sampling scheme. *Ann. Statist.* **21** 1567–1590. MR1241279

Zhang, H. (2004). Inconsistent estimation and asymptotically equal interpolations in model-based geostatistics. *J. Amer. Statist. Assoc.* **99** 250–261. MR2054303

Zhang, H. and Zimmerman, D. L. (2005). Towards reconciling two asymptotic frameworks in spatial statistics. *Biometrika* **92** 921–936. MR2234195

J. Du
V. S. Mandrekar
Department of Statistics and Probability
Michigan State University
East Lansing, Michigan 48824
USA
E-mail: dujuan@stt.msu.edu
    mandrekar@stt.msu.edu

H. Zhang
Department of Statistics
Purdue University
West Lafayette, Indiana 47907
USA
E-mail: zhanghao@purdue.edu